# The butterfly sequence: the second difference sequence of the numbers of integer partitions with distinct parts, its pentagonal number structure, its combinatorial identities and the cyclotomic polynomials $1-x$ and $1+x+x^2$

*Cristiano Husu*

*Department of Mathematics*

*University of Connecticut*

*UConn, Stamford, CT 06901*

*To B. and W.*

## Abstract

Based on the author's previous work on the Jacobi identity for relative twisted vertex operator algebras and modules, and on the generating function identities for affine Lie algebras, we interpret the second difference of the sequence of the number of integer partitions with distinct parts (the strict partitions), for $n \geq 5$, as the sequence of the strict partitions of $n$ with *at least three parts, the three largest parts consecutive*, and *the smallest part at least two*. The name *butterfly* describes both the sequence's interpretation and the underlying bijection between the set of strict partitions of each positive integer $m$, $m \geq 5$, with the two largest parts consecutive, and a subset of the same kind of strict partitions of $m+1$. Using the cyclotomic polynomials $1-x$ and $1+x+x^2$, we compute generating function identities of the butterfly sequence and related sequences both as infinite products and as series filtered by the number of parts of the corresponding partitions, and we see that the number of partitions of positive integers $n$ with odd parts greater than or equal to 5 is the sum of three consecutive terms of the butterfly sequence. We also determine a subtler *merging and splitting* construction of the butterfly sequence as a sequence of some of the partitions with odd parts larger or equal to 3, and we offer a related detailed interpretation of the butterfly sequence as a sequence of what we define as *generalized pentagonal*, *pentagonal "with domino"*, and *non-pentagonal butterfly partitions*. Finally, Euler's Pentagonal Number Theorem and a slightly different specialization of the Jacobi Triple Product lead to recursive algorithms to compute the butterfly sequence and related sequences using *pentagonal number* sequences and *the series of triangular powers*.

## INTRODUCTION

The present article is an extended version of the paper titled "the butterfly sequence: the second difference sequence of the numbers of integer partitions with distinct parts and its pentagonal number structure", [H], a study of integer partitions, a mix of elementary and subtler ideas, based on the author's previous work on the Jacobi identity for relative twisted vertex operator algebras and modules, and on the generating function identities for affine Lie algebras. With respect to the compact version [H], the present extended version contains further basic details of the notation, definitions, propositions, theorems and proofs which are also found in [H]. It also contains several introductory numerical examples that are not in [H], and eleven additional figures (the Young diagram figures 1 - 5 at the end of the present introduction, and the Young diagram figures 7 – 13 in sections 1 and 4, here) not in [H]. To help readers of the compact version [H] find the corresponding extended details in the present extended version, we have included the present paper, in red, the corresponding numbering in [H] of definitions, propositions, theorems and proofs.

Most of the work, in both [H] and the present extended version, is about what we call the *butterfly sequence.* This is defined as the second difference of the sequence of the number of integer partitions with distinct parts (the second difference of the sequence of the number of strict partitions). The name of the sequence is for the interpretation of the sequence as the (first) difference sequence of the strict partitions with the two largest parts consecutive: in the Young diagrams associated with this interpretation, two squares joined at a vertex form a rough sketch of a butterfly. The result of our interpretation is that the terms of the butterfly sequence, for $n \geq 6$, count the number of strict partitions of $n$ with at least three parts, the three largest parts consecutive, and the smallest part at least 2 (the *butterfly partitions,* in short).

The second difference sequence of the integer partitions with distinct parts came to the attention of the author of the present paper when comparing two very specific combinatorial structures. The first is the structure of the level 2 standard representations of the affine Lie algebra $A_1^{(1)}$ (see the last two sections of [LW1]), which is directly related to the butterfly sequence as the Weyl character formula, for these representations, gives Euler's identities for integer partitions with distinct even parts (the 'domino version" of the integer partitions with distinct parts, here), and their counterparts for distinct odd parts. The second combinatorial structure, more complicated and harder to visualize because of the intertwined vertex operators, is that of the Jacobi identity for relative twisted vertex operators, and the combinatorial shifts of the Jacobi identity determined by multiplication (and division) of the cyclotomic polynomials in the square root of the formal variable for the affine Lie algebra $A_1^{(1)}$.

These combinatorial structures of the affine Lie algebra $A_1^{(1)}$ are, in turn, based on the more general work on vertex operator algebras and modules, and on the general representations of affine Lie algebras. This work includes the construction of C. Y. Dong and J. Lepowsky, [DL1] and [DL2], of the Jacobi identity for relative vertex operators and for relative twisted vertex operators, the author's previous work, [H1] and [H2], on the Jacobi identity for relative twisted vertex operators and on the generating function identities for standard representations of the affine Lie algebras $A_1^{(1)}$ and $A_2^{(2)}$, and the general Z-operator construction of J. Lepowsky and R. Wilson, [LW1] and [LW2], of the standard modules of affine Lie algebras, Euler (distinct even and distinct odd parts) partition identities and the Rogers-Ramanujan identities. In all these vast areas of mathematics, a *striking persistence* of similar combinatorial computations occurs in three algebra structures. The first structure, in order of generality, is that of the vertex operator algebras and modules in which the Jacobi identity is mainly manipulated via the *main algebraic property of the Dirac delta function*

$$f(x)\delta(x) = f(1)\delta(x)$$

for $f(x)$ in $\mathbb{C}[[x]]$, and $\delta(x) = \sum_{m=-\infty}^{\infty} x^m = \frac{1}{1-x} + \frac{\frac{1}{x}}{1-\frac{1}{x}}$ (in the author's previous work, the formal variable x is $\frac{z_2^{1/2}}{z_1^{1/2}}$, for $A_1^{(1)}$, and $\frac{z_6^{1/2}}{z_6^{1/2}}$, for $A_2^{(2)}$). The second structure is that of the affine Lie algebras and modules. For $A_1^{(1)}$ and $A_2^{(2)}$, *commutators and anticommutators of the algebra generating functions* are computed from the Jacobi identity by extracting the residue with respect to the *binding formal variable* and by further multiplication (and division) of the Jacobi identity by the cyclotomic polynomials $\Phi_1(x) = 1 - x$, $\Phi_2(x) = 1 + x + x^2$, $\Phi_3(x) = 1 + x = \Phi_1(-x)$, and $\Phi_6(x) = 1 + x + x^2 = \Phi_2(-x)$.

The third structure is that of the partitions and partition identities: *Euler's distinct even parts and distinct odd parts identities* for the principal character of level 2 standard $A_1^{(1)}$-modules; the *Rogers-Ramanujan identities* for level 3 standard $A_1^{(1)}$-modules; the *domino Rogers-Ramanujan identities* for level 2 standard $A_2^{(2)}$-modules ("domino" for the Young diagram domino visualization in place of each square); the *Capparelli identities* (see [C1], [C2] and [TX]) for level 3 standard $A_2^{(2)}$-modules; and further partitions identities to be determined (see [BM] for standard $A_7^{(2)}$-modules, and [B] for partitions of the principal character of standard modules of the affine Lie algebras $A_7^{(2)}$, $C_3^{(1)}$, $F_4^{(1)}$, and $G_2^{(1)}$).

As for a detailed description of the content of the present paper, in the first section, we recall the sequence of strict partitions, $q(n)$, its (first) difference sequence, $r(n)$, and several interpretations of the difference sequence. We, then, define the *butterfly sequence*, $s(n)$, for $n \geq 0$, as the second difference sequence of the number of strict partitions, and we give a detailed description of the sets of partitions and the bijections between these sets that lead, for $n \geq 6$, to the main interpretation of the butterfly sequence (Proposition 1.3, and the related corollaries 1.4, 1.5 and 1.6, here). In the second section, we construct generating functions of the butterfly sequence and related sequences both as infinite products and as series filtered by the number of parts of the corresponding partitions. The generating function identities, in turn, display the relationship between the butterfly sequence and the partitions of $n$ with odd parts larger or equal to 5 determined by the cyclotomic polynomial $1 + x + x^2$, and set a partition identity between the number of partitions with odd parts larger or equal to 5, and the number of partitions with distinct parts, at least three parts larger or equal to 2, up to two 1 part (the part 1 may not occur, or occur once, or occur twice), and the three largest parts consecutive (Proposition 2.1, here). The third section, here, is more complicated, as it contains a subtler *merging and splitting* construction of the butterfly sequence as a sequence of partitions with odd parts larger or equal to 3 (Propositions 3.4 and 3.5, here, and, respectively, Propositions 3.3 and 3.4 in [H]). This construction requires the further distinction of the partitions of the butterfly sequence in terms of the partitions with the second largest part even, and their odd counterparts. The fourth section is related to the third, and contains a detailed interpretation of the partitions of the butterfly sequence in terms of *generalized pentagonal numbers*, *generalized pentagonal numbers with two extra units*, and what we define as *generalized pentagonal, generalized pentagonal with domino*, and *non-pentagonal butterfly partitions* (Theorem 4.6, and the related corollaries 4.6, 4.7, 4.8 and 4.9, here). In the fifth, and final, section, we use Euler's Pentagonal Number Theorem, and a second specialization of the Jacobi Triple Product, to construct the terms of the butterfly sequence and of the related sequences as pentagonal and triangular number sequences of the partition function. We conclude the fourth section with recursive algorithms to compute the butterfly sequence and the related sequences, algorithms that use only pentagonal sequences (in the sequence inputs) and triangular numbers (in the sequence output). In these algorithms, we see how the multiplication of the series of triangular powers (the series on the right side of our second specialization of the Jacobi Triple Product) by the cyclotomic polynomials $1 - x$ and $1 + x + x^2$ folds and scrambles the first few terms of the butterfly sequence and of the sequence of partitions with odd parts larger or equal to 5.

We conclude the introduction with a note about the first insight that led the present paper. In view of previous work on vertex operators, affine Lie algebras and combinatorial identities, *integers are interpreted in different ways within the same mathematical structure*. In the present paper, this general principle resonates as positive integers $n$ are interpreted in terms of a specific subset of numbers (here, *extended pentagonal numbers* as inputs of sequences of strict partitions), and in terms of specific partition structures of that the numbers embody (here, *butterfly generalized pentagonal partitions, butterfly pentagonal partitions with dominoes*, and *butterfly non-pentagonal partitions* as a specific interpretation of the outputs of the butterfly sequence,). A few examples, together with the corresponding Young diagram should clarify this insight (for the related general pentagonal number structure, see section 4, below).

**Example 1:** $n = 9$ is the smallest integer for which a nonzero butterfly partition exists, and indeed 9 is the *extended pentagonal number* $7 = 5 + 2$, *with the additional* 2 *for the second difference structure*, and, at the same time, 9 is interpreted as the parts in the (butterfly) partition $4 > 3 > 2$.

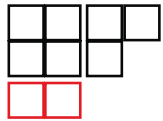

**Figure 1:** the butterfly partition $4 > 3 > 2$. In the diagram, the extended pentagonal number, 7, is presented as the black 2x2 square together with the second triangular number on the right. The additional red domino is generated by the second difference sequence.

**Example 2:** $n = 12$ is the pentagonal number $\frac{3}{2} \, 3^2 - \frac{1}{2} \, 3$ which is also interpreted as the result of the butterfly partition $5 > 4 > 3$.

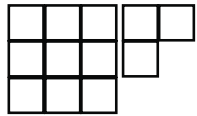

**Figure 2:** the butterfly partition $5 > 4 > 3$. The pentagonal number 12 is presented as a 3x3 square together with the second triangular number.

**Example 3:** $n = 14$ is the pentagonal number 12 with the additional 2, and is also the sum of the parts in the partition $5 > 4 > 3 > 2$.

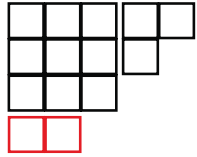

**Example 4:** $n = 15$ is the extended pentagonal number $\frac{3}{2} \, 3^2 + \frac{1}{2} \, 3$, and the sum of the parts in the partition $6 > 5 > 4$. (We think of $6 > 5 > 4$ as the partition in Example 2, extended by the addition of a unit to each part.

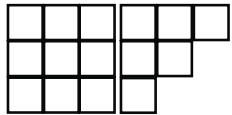

**Example 5:** $n = 17$ is the extended pentagonal number 15 with the additional 2, and result of the partition $6 + 5 + 4 + 2$.

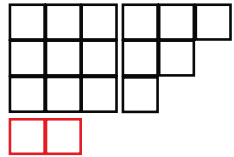

**Example 6:** $n = 18$ is the smallest positive integer that has two distinct (nonzero) butterfly partitions. We think of $18 = 15 + 3 = 12 + 3 + 3$ is as a *twice-extended pentagonal number* (for the repeated addition of three units to the pentagonal number 12). The two butterfly partitions $7 + 6 + 5$, and $6 + 5 + 4 + 3$ are associated to each other by their relationship with the partition in Example 4: the first is obtained by adding one unit to each part of $6 + 5 + 4$, the second by adding the part 3. Note that in the first partition, the second part is even, while, in the second partition, the second part is odd.

(Examples $1 - 5$, and the corresponding Young diagrams, do not appear in [H]. On the other hand, Example 6, above, also appears at the end of the introduction in [H], p. 2, and Figure 6, below, is also Figure 1 in [H], p. 2.)

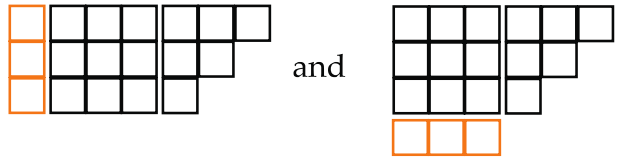

**Figure 6:** the two butterfly partitions 7 > 6 > 5 and 6 > 5 > 4 > 3. The two partitions are associated with each other by the orange bars (vertical, on the left, and horizontal, on the right) around the identical pentagonal black 6+5+4 core. The reader should imagine that the partition on the right is created from the partition on the left, as a unit is removed from each of the (consecutive) parts and a 3 part is inserted in place of the three units, and, vice versa, the partition on the left is created from the partition on the right, as the horizontal part 3 is removed and a unit is added to each of the remaining parts.

# 1 THE BUTTERFLY SEQUENCE, THE SECOND DIFFERENCE SEQUENCE OF THE NUMBER OF INTEGER PARTITIONS WITH DISTINCT PARTS

Let $\{q(n)\}$, for $n$ nonnegative integers, be the sequence of the number of *strict partitions* of $n$, which is to say $q(0) = 1$, and, for $n \geq 1$, $q(n)$ is the sequence of the number of partitions of $n$ with distinct parts. The sequence $\{q(n)\}$ can be interpreted in many other ways, for some of which the reader can quickly check the sequence A000009 in the On-Line Encyclopedia of Integer Sequences (O. E. I. S.). Here are the first few terms of the sequence for some insights about the background of the present paper:

$$\{q(n)\} = \{1, 1, 1, 2, 2, 3, 4, 5, 6, 8, 10, 12, 15, 18, 22, 27, 32, 38, 46, 54, 64, 76, 89, 104 \dots \}, \quad (1.1)$$

for $n = 0, 1, 2, 3, \dots$.

For instance, we note the *(generalized) pentagonal number structure* of the inputs corresponding to the odd outputs in (1.1), a straightforward consequence of Euler's insight of dividing the strict partitions into two sets, the set of partitions with an even number of distinct parts, and the set of partitions with an odd number of distinct parts: for $n \geq 1$, so that $q(n)$ is even *except when $n$ is a generalized pentagonal number* (cf. [A], Theorem 1.6 and Corollaries 1.7 and 1.8).

Let $\{r(n)\}$, for $n \geq 0$, be the difference sequence of $\{q(n)\}$: set $q(-1) = 0$ and define

$$\{r(n)\} = \{dq(n)\} = \{q(n) - q(n-1)\} \quad (1.2)$$

for $n \geq 0$. Here are the first few terms of the difference sequence for a quick introduction to the main insight of the present section:

$$\{r(n)\} = \{\,1, 0, 0, 1, 0, 1, 1, 1, 1, 2, 2, 2, 3, 3, 4, 5, 5, 6, 8, 8, 10, 12, 13, 15, 18, 20, 23, 27, 30,$$
$$34, 40, 44, 50, 58, 64, 73, 83, 92, 104, 118, 131, \dots \,\}$$

for $n = 0, 1, 2, 3, \dots .$

As for (1.1) the sequence (1.2) can be interpreted in many ways, several of which in the O. E. I. S., sequence A087897. For the purposes of the present work, we mention three interpretations of the difference sequence (1.2) closely related to the above *strict partitions* definition of (1.1). The first interpretation is that, for $n \geq 3$, $r(n)$ is the sequence of the *strict partitions* of $n$ *with the two largest parts consecutive*, that is to say, $r(n)$ is the sequence of the partitions of $n$ with distinct parts, and the two largest parts consecutive. The second is that, for $n \geq 3$, $r(n)$ is the sequence of the *partitions of $n$ into odd parts greater than 3*. And the third is that, for $n \geq 4$, $r(n)$ is the sequence of *partitions of $n$ with distinct parts which are not powers of two.*

As for the pentagonal number structure of (1.1), the definition of (1.2) gives what we call the *first difference pentagonal number structure* of (1.2): for $n \geq 3$, $r(n)$ is even *except when $n$ is a generalized pentagonal number, or a generalized pentagonal number plus one.*

The interpretation of the difference sequence (1.2), for $n \geq 3$, as the sequence of partitions with distinct parts and the two largest parts consecutive can be easily seen by means of the bijection (in Figure 7 Young diagrams, below) between the set A of all the strict partitions of $n - 1$, for $n \geq 4$, and the set B of the strict partitions of $n$, for $n \geq 4$, with the two largest parts of difference greater or equal to two units, obtained by adding 1 unit to the largest part of every element of A to create a partition of B, and, vice versa, removing a unit from each partition of B to create a partition of A.

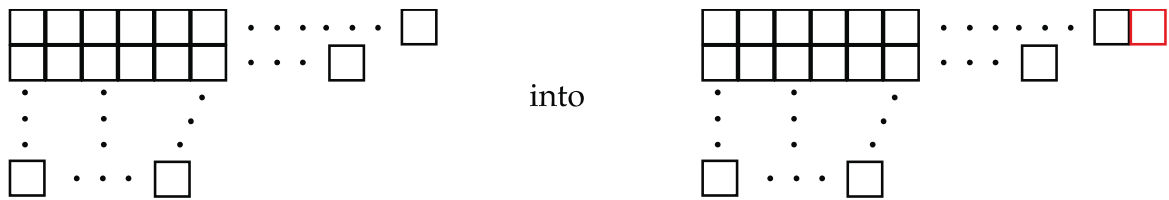

**Figure 7:** the bijection that gives the interpretation of (1.2) as the sequence of the number of strict partitions of $n$, for $n \geq 3$, with the two largest parts consecutive. The left diagram represents an arbitrary strict partition of $n - 1$, for $n \geq 4$, and the diagram on the right, with the additional red unit added to the largest part, represents the result of the bijection, a strict partition of $n$ with the two largest parts of difference greater than or equal to 2.

With the interpretation of (1.2) as the sequence of strict partitions with the two largest parts consecutive, we divide the set of the strict partitions with the two largest parts consecutive into the set $R_1$ of such partitions with the smallest part at least 2, and the set $R_2$ of such partitions with the smallest part exactly 1, and we consider the corresponding sequences. For $n \geq 3$, we define the sequence $r_1(n)$ of the number of strict partitions of $n$ with the two largest parts consecutive and the smallest part at least 2. Here are the first few terms of the sequence (note that we start the start the sequence at $n = 3$ although the first nonzero element of $R_1$, $5 = 3 + 2$, occurs for $n = 5$):

$$\{r_1(n)\} = \{\ 0, 0, 1, 0, 1, 0, 2, 0, 2, 1, 2, 2, 3, 2, 4, 4, 4, 6, 6, 7, 8, 10, 10, 13, 14, 16, 18, \dots\} \quad (1.3)$$

for $n = 3, 4, 5, 6, 7, \dots$, and, again for $n \geq 3$, we define the sequence $r_2(n)$ of the number of strict partitions of $n$ with the two largest parts consecutive and the smallest part exactly 1:

$$\{r_2(n)\} = \{\ 1, 0, 0, 1, 0, 1, 0, 2, 0, 2, 1, 2, 2, 3, 2, 4, 4, 4, 6, 6, 7, 8, 10, 10, 13, 14, 16, \dots\} \quad (1.4)$$

for $n = 3, 4, 5, 6, 7, \dots$. Obviously,

$$r_1(n) + r_2(n) = r(n), \quad for\ n \geq 3, \quad (1.5)$$

and also note that

$$r_1(n-1) = r_2(n), \quad for\ n \geq 4, \quad (1.6)$$

as every partition of $n - 1$ in $R_1$ corresponds to a partition of $n$ in $R_2$ by adding the part 1, and, vice versa, by removing the part 1 from every partition of $n$ in $R_2$ to reconstitute a partition of $n - 1$ in $R_1$.

To introduce a slightly subtler bijection, we also define the following two sets $R_1'$ and $R_1''$. The first, $R_1'$, is the set of the strict partitions *with the two largest parts consecutive, the smallest part at least 2,* and with the further condition that, *for partitions with at least three parts, the second largest part is at least two units greater than the third largest part.* The second, $R_1''$, is the set of the strict partitions *with at least three parts, the three largest parts consecutive, and the smallest part at least 2.* We also consider the corresponding sequences of the number of partitions of $n$, for $n \geq 5$, in $R_1'$ and $R_1''$, respectively,

$$\{r_1'(n)\} = \{1, 0, 1, 0, 2, 0, 2, 1, 2, 2, 3, 2, 4, 4, 2, 6, 6, 7, 8, 10, 10, 13, 14, \dots\}, \quad n = 5, 6, 7, \dots, \quad (1.7)$$

and, respectively,

$$\{r_1''(n)\} = \{\ 0, 0, 0, 0, 1, 0, 0, 1, 0, 1, 1, 0, 1, 2, 0, 2, 2, 1, 2, 3, 2, 3, 4, \dots\}, \quad n = 5, 6, 7, \dots. \quad (1.8)$$

Note that, for $n \geq 5$,

$$r_1'(n) + r_1''(n) = r_1(n). \quad (1.9)$$

Also note that we choose to start both sequences at $n = 5$, although the first nonzero element of $R_1''$, the partition $9 = 4 + 3 + 2$, occurs for $n = 9$.

We finally introduce the second difference sequence of the sequence of strict partitions.

**Definition 1.1:** Let $\{s(n)\}$, for any nonnegative integer $n$, be the difference sequence of $\{r(n)\}$ as in (1.2), or, equivalently, the second difference sequence of $q(n)$ as in (1.1): set $r(-1) = q(-1) = q(-2) = 0$, and define, for $n \geq 0$,

$$\{s(n)\} = \{dr(n)\} = \{r(n) - r(n-1)\}, \tag{1.10}$$

or, equivalently,

$$\{s(n)\} = \{d^2 q(n)\} = \{q(n) - 2q(n-1) + q(n-2)\}, \tag{1.11}$$

also for $n \geq 0$.

(Definition 1.1, here, is also <span style="color:red">Definition 1.1</span> in <span style="color:red">[H]</span>, p. 3.)

□

Here are the first fifty-two terms of the sequence (1.10) to introduce the reader to the pentagonal structure of the sequence:

$$\{s(n)\} = \{\ 1, -1, 0, 1, -1, 1, 0, 0, 0, 1, 0, 0, 1, 0, 1, 1, 0, 1, 2, 0, 2, 2, 1, 2, 3, 2, 3, 4, 3, 4, 6, 4, 6,$$
$$8, 6, 9, 10, 9, 12, 14, 13, 14, 19, 18, 22, 26, 24, 30, 34, 34, 40, 45, \dots \},$$

for $n = 0, 1, 2, 3, 4, \dots$ .

As a consequence of the pentagonal number structure in (1.1), and as for what we call the first difference pentagonal number structure in (1.2), we can make a similar observation, in Remark 1.2, below, about the second difference pentagonal number structure in the sequence (1.10).

**Remark 1.2:** For $n \geq 5$ and $n \neq 7$, note that $s(n)$ is even *except when $n$ is a generalized pentagonal number, or a generalized pentagonal number plus two (or a generalized pentagonal number "with domino", that is $s(n)$ is even except when*

$$n = \ k^2 + \tfrac{1}{2}k(k-1) = \tfrac{3}{2}\,k^2 - \tfrac{1}{2}k, \text{ with } k = \pm 3, \pm 4, \pm 5, \dots ,$$

or

$$n = \ k^2 + \tfrac{1}{2}k(k-1) + 2 = \tfrac{3}{2}\,k^2 - \tfrac{1}{2}k + 2, \text{ with } k = -2, \pm 3, \pm 4, \pm 5, \dots .$$

We will refer to this feature of the sequence (1.2) as the *second difference generalized pentagonal number structure,* or just the *generalized pentagonal number structure,* of the sequence (1.10). (The exception of the first eight terms of the sequence (1.10), the terms $s(n)$, for $n \leq 7$, is the result of a sort of *overlap,* or *wrinkling,* of the corresponding term in (1.1) and (1.2). The algebraic interpretation of this exception is at the end of section 2, in the identities (2.11) and (2.12), here.)

(Remark 1.2, here, is also <span style="color:red">Remark 1.2</span> in <span style="color:red">[H]</span>, p. 3.)

□

The following proposition offers the main interpretation of the sequence (1.10). The proof of the proposition, below, contains the bijection (Figure 8, below) for which the term "butterfly" in the present work is derived.

**Proposition 1.3:** *For $n \geq 6$, each term of the sequence (1.10), is equal to the corresponding term of the sequence (1.8).*

*In other words, the sequence (1.10), for $n \geq 6$, can be interpreted as the sequence of the number of partitions of $n$ with distinct parts $p_i, i = 1, 2, \ldots, k$, such that $k \geq 3$, and*

$$p_1 = p_2 + 1 = p_3 + 2 > p_4 > \cdots > p_k > 1. \qquad (1.12)$$

*Proof.* For $n \geq 6$, using the Definition 1.1 together with the identities (1.5), (1.6) and (1.9), we see that

$$s(n) = r(n) - r(n-1) =$$
$$= r_1(n) + r_2(n) - r_1(n-1) - r_2(n-1) = \qquad (1.13)$$
$$= r_1(n) - r_2(n-1) = r_1'(n) + r_1''(n) - r_2(n-1).$$

Moreover, for $n \geq 6$, we can and do inject the partitions of $n-1$ in the set $R_2$ into $R_1'$, by removing the smallest part (the part equal to 1) in each partition of $n-1$ in $R_2$, and, to create a partition of $n$ in $R_1'$, by adding one unit to each of the two largest parts of the result (see the Young diagrams in Figure 8, below), and, vice versa, we inject the partitions of $n$ in the set $R_1'$ into $R_2$, by removing one unit from each of the two largest parts of each partition of $n$ in $R_1'$, and, to create a partition of $n-1$ in $R_2$, adding a part equal to 1 to the result.

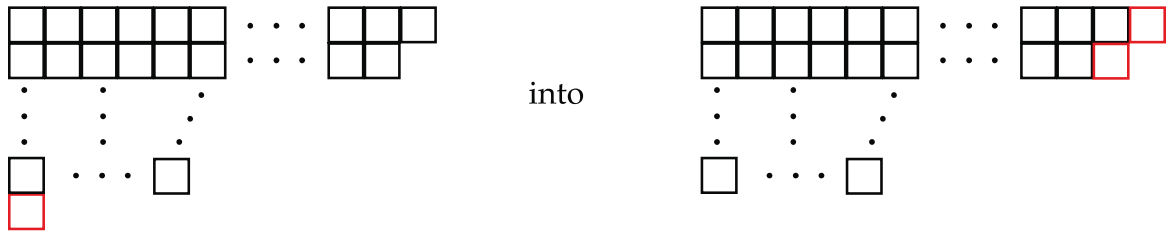

into

**Figure 8:** the *butterfly bijection*: the bijection between the set of strict partitions of $n-1$, for $n \geq 6$, with at least three parts, the two largest parts consecutive, and the smallest part 1, and the set of strict partitions of $n$, for $n \geq 6$, with at least three parts, the two largest parts consecutive, the third largest part at least two units smaller than the second largest part, and the smallest part greater than 1. The above diagram yields the interpretation of (1.10) as the sequence of the number of strict partitions of $n$, for $n \geq 6$, with at least three parts, the three largest parts consecutive, and the smallest part at least 2. The Young diagram on the left represents a strict partition of $n-1$, for $n \geq 6$, with at least three parts, the two largest parts consecutive, and the smallest part 1, while the diagram on the right represents a strict partition of $n$, for $n \geq 6$, with at least three parts, the two largest parts consecutive, the third largest part at least two units smaller than the second largest part, and the smallest part greater than 1. The bijection is obtained by removing the part 1 from the partition on the left, and, to create the partition on the right, adding a unit to each of the largest parts, and, vice versa, removing a unit from each of the two largest parts on the left, and, to create the partition on the right, adding the part 1.

Through the above bijection, we see that, for $n \geq 6$,

$$r_2(n-1) = r_1'(n),$$

and, as for (1.13),

$$s(n) = r_1''(n).$$

□

(Proposition 1.3, here, is also <span style="color:red">Proposition 1.3</span> in <span style="color:red">[H]</span>, p. 3.)

We name the bijection in the proof of Proposition 1.3 and Figure 8, above, the *butterfly bijection* for the red squares (a sketch of a butterfly) in the Young diagrams of Figure 8 (or, if the reader wishes, a square that turns into a rough sketch of a butterfly). By the way, similar "butterfly" diagrams can be used to visualize Euler's pentagonal number theorem (cf. [A2], Theorem 1.6 and Corollaries 1.7 and 1.8). For this bijection, we will also refer to the sequence (1.10) as the *butterfly sequence,* and, for $n \geq 9$, we will also refer to the partitions (1.12) in Proposition 1.1 of the as the *butterfly partitions ($n \geq 9$ as the first nonzero such partition is $9 = 4 + 3 + 2$).*

Note that the *trapezoidal partition structure*

$$p_1 = p_2 + 1 = p_3 + 2,$$

for $p_3 = 2, 3, 4, \ldots$, of the three largest parts of the partitions in Proposition 1.1, can obviously be replaced with the rectangular structure of three equal parts,

$$p_1 = p_2 = p_3,$$

for $p_3 = 3, 4, 5, \ldots$, by removing one unit from the largest part, and adding one unit to the third largest part, thus obtaining the following corollary.

**Corollary 1.4:** *The sequence (1.10), for $n \geq 6$, can be interpreted as the sequence of the number of partitions of $n$ with parts $p_i$, $i = 1, 2, \ldots, k$, such that (a) $k \geq 3$, (b) the first three largest parts are equal, (c) the remaining parts, if present, are distinct, (d) the fourth largest part, if present, is two units smaller than the third largest part, and (e) the smallest part is not 1.*

(Corollary 1.4, here, is also <span style="color:red">Corollary 1.4</span> in <span style="color:red">[H]</span>, p. 3.)

☐

Switching the number of parts with the sizes of the parts in the butterfly partitions (1.12) of Proposition 1.3, we have the following corollary.

**Corollary 1.5:** *The sequence (1.10), for $n \geq 6$, can be interpreted as the sequence of the number of partitions of $n$ with parts $q_j$, $j = 1, 2, \ldots, h$, such that (a) $h \geq 4$, (b) the parts $q_j$, $j = 1, 2, \ldots, h$, can be arranged in non-increasing order so that each pairs of consecutive parts are either equal or differ by one unit, (c) the first two largest parts are equal, and (d) the three smallest parts are 3, 2, and 1. In detailed notation, the sequence (1.10), for $n \geq 6$, can be interpreted as the sequence of the number of partitions of $n$ with parts $q_j$, $j = 1, 2, \ldots, h$, of the form*

$$q_1 = q_2 \geq q_3 \geq q_4 \geq \cdots \geq q_{h-2} = 3 > q_{h-1} = 2 > q_h = 1, \tag{1.14}$$

*with $h \geq 4$, $q_j - q_{j+1} \leq 1$, for all $j = 3, 4, \ldots, h - 3$.*

(Corollary 1.5, here, is also <span style="color:red">Corollary 1.5</span> in <span style="color:red">[H]</span>, p. 3.)

☐

And, similarly, switching the number of parts with the sizes of the parts in the partitions of Corollary 1.4, we have the following corollary.

**Corollary 1.6:** *The sequence (1.10), for $n \geq 6$, can be interpreted as the sequence of the number of partitions of $n$ with parts $q_j$, $j = 1, 2, \ldots, h$, such that (a) $h \geq 3$, (b) the parts $q_j$, $j = 1, 2, \ldots, h$, can be arranged in non-increasing order so that each pairs of consecutive parts are either equal or differ by one unit, (c) the first two largest parts are equal, and (d) the two smallest* *parts are both equal to 3, 2. In detailed notation, the sequence (1.10), for $n \geq 6$, can be interpreted as the sequence of the number of partitions of $n$ with parts $q_j$, $j = 1, 2, \ldots, h$, of the form*

$$q_1 = q_2 \geq q_3 \geq q_4 \geq \cdots \geq q_{h-1} = q_h = 3, \tag{1.15}$$

*with $h \geq 3$, $q_j - q_{j+1} \leq 1$, for all $j = 3, 4, \ldots, h-2$.*

(Corollary 1.6, here, is also <span style="color:red">Corollary 1.6</span> in <span style="color:red">[H]</span>, p. 4.)

$\square$

We conclude the section with a note about three subsequences of (1.10):

$$\{s(3m)\} = \{\ 1, 1, 1, 2, 2, 3, 4, 6, 8, 10, 14, 19, 26, 34, 45, \ldots \},$$

$$\{s(3m+1)\} = \{\ 0, 0, 0, 0, 1, 2, 3, 4, 6, 9, 13, 18, 24, 34, \ldots \},$$

and

$$\{s(3m+2)\} = \{\ 0, 1, 1, 2, 2, 3, 4, 6, 9, 12, 14, 22, 30, 40, \ldots \},$$

for $m = 3, 4, 5, \ldots$ . These subsequences are obviously also part of the interpretation provided by Proposition 1.3. Their butterfly partitions, or, to be more specific, the three largest parts of their butterfly partitions, account for the non-decreasing pattern of the three subsequences which is also easily notice in the first few terms.

## 2 GENERATING FUNCTIONS, INFINITE PRODUCTS, FILTRATION SERIES, AND THE RELATIONSHIP BETWEEN THE BUTTERFLY SEQUENCE THE SEQUENCE OF PARTITIONS WITH ODD PARTS LARGER OR EQUAL TO 5

In the second section, and in the rest of this paper, we use the formal variable $x$, in place of the most commonly used variable $q$, to write down generating functions of partition sequences. From the author's point of view, this is both to be consistent with related previous work notation in [H1] and [LW2], and to avoid confusion with $q(n)$, the number of strict partitions of $n$.

We begin by recalling the generating function of $q(n)$:

$$\prod_{n=1}^{\infty} (1 + x^n) = 1 + \sum_{k=1}^{\infty} \frac{x^{\frac{k(k+1)}{2}}}{(1-x)(1-x^2)\ldots(1-x^k)}, \qquad (2.1)$$

which, using the difference of squares, is apparently also the generating function of partitions with odd parts (hence showing that, for any $n$, the number of partitions with odd parts is equal to the number of strict partitions):

$$\prod_{h=0}^{\infty} \frac{1}{1-x^{2h+1}} = 1 + \sum_{k=1}^{\infty} \frac{x^{\frac{k(k+1)}{2}}}{(1-x)(1-x^2)\ldots(1-x^k)}. \qquad (2.2)$$

Note that each term of the infinite series in (2.1) and (2.2), for each $k > 0$, represents the strict partitions with exactly $k$ number of parts, so that the series on the right sides of (2.1) and (2.2), respectively, are the *filtrations* of the strict partitions and the strict partitions with the two largest parts consecutive, filtered according to the number of parts. The first term, 1, representing the number of (trivial) partition of zero, is what we see as a "ground", or "ghost" partition.

In a similar way, the generating function of the number of strict partitions with the largest two parts consecutive (see (1.2)) can be written as

$$(1-x)\prod_{n=1}^{\infty} (1 + x^n) = 1 + \sum_{k=2}^{\infty} \frac{x^{\frac{k(k+1)}{2}}}{(1-x^2)\ldots(1-x^k)}, \qquad (2.3)$$

in which the left side, again using the difference of squares, can also written as the generating function of partitions with odd parts *larger than* 1 (hence showing that, for any n, the number of partitions with odd parts larger than 1 is equal to the number of strict partitions with consecutive largest parts)

$$\prod_{h=1}^{\infty} \frac{1}{1-x^{2h+1}} = 1 + \sum_{k=2}^{\infty} \frac{x^{\frac{k(k+1)}{2}}}{(1-x^2)\ldots(1-x^k)}. \qquad (2.4)$$

Again, each terms of the infinite series on the right side of (2.3), for each $k > 1$, represents the strict partitions with the two largest parts consecutive and exactly $k$ number of parts, so that the series in (2.4) is the *filtration* of the strict partitions with the two largest parts consecutive, filtered according to the number of parts. The first term, 1, again, represents the number of (trivial) partition of zero. Note that the multiplication by (1 - $x$) of the first two terms of the right side of (2.1) and (2.2) leaves the zero coefficient of $x$ as for there are not one part partitions with consecutive parts.

For the butterfly sequence, further multiplication of (2.3) by the factor $(1 - x)$ leads to the second difference counterparts of (2.1) and (2.2):

$$(1-x)^2 \prod_{n=1}^{\infty} (1+x^n) = (1-x)^2 + (1-x)^2 \sum_{k=1}^{\infty} \frac{x^{\frac{k(k+1)}{2}}}{(1-x)(1-x^2)\dots(1-x^k)}, \qquad (2.5)$$

and

$$(1-x) \prod_{h=1}^{\infty} \frac{1}{1-x^{2h+1}} == 1 - x + \frac{1}{(1+x)} \sum_{k=3}^{\infty} \frac{x^{\frac{k(k+1)}{2}}}{(1-x^3)(1-x^4)\dots(1-x^k)}. \qquad (2.6)$$

Note that the triangular exponents of $x$ in the series of (2.5) and (2.6) give an algebraic insight in the triplet structure of the terms of the sequence (1.10) (see the note at the end of section 1, here).

The infinite product on the left side can also be written as

$$\frac{1}{(1+x+x^2)} \prod_{h=2}^{\infty} \frac{1}{1-x^{2h+1}}, \qquad (2.7)$$

and, as such, uncovers the relationship between the butterfly sequence and the number of partitions with odd parts larger or equal to 5. The infinite series, the right side of (2.5) and (2.6), requires algebraic manipulations to display the *filtration* of the butterfly sequence in terms of equal numbers of parts. The algebraic manipulations that the author found are complicated. It is easier, as shown in the rest of this section, to compute the filtration directly for $n \geq 9$, and, for the first nine terms, to take a quick look at the partitions with odd parts larger or equal to 5.

The butterfly sequence partitions of the form $4 + 3 + 2, 5 + 4 + 3, 6 + 5 + 4, \dots$, for $n \geq 9$, give the terms of the generating function $x^9 + x^{12} + x^{15} + \cdots = \frac{x^9}{1-x^3}$, and, similarly, for four parts, $\frac{x^{14}}{1-x^3} + \frac{x^{18}}{1-x^3} + \frac{x^{22}}{1-x^3} + \cdots = \frac{x^{14}}{(1-x^3)(1-x^4)}$, and so forth, so that the generating function of the butterfly sequence can be written, *for $n \geq 9$*, as the filtered series

$$\sum_{k=3}^{\infty} \frac{x^{\frac{k(k+3)}{2}}}{(1-x^3)(1-x^4)\dots(1-x^k)}. \qquad (2.8)$$

For the first nine terms of the butterfly sequence, consider the generating function of the sequence of partitions of $n$, *for $n \geq 9$*, with odd parts larger or equal to 5. The infinite product side, (2.7), multiplied by the polynomial $1 + x + x^2$, and the series (2.8) also multiplied by the polynomial $1 + x + x^2$, give the following identity: *the coefficient of $x^n$, for $n \geq 9$, in*

$$\prod_{h=2}^{\infty} \frac{1}{1 - x^{2h+1}}$$

is equal to *the coefficient of $x^n$, for $n \geq 9$, in*

$$\sum_{k=3}^{\infty} \frac{(1 + x + x^2) x^{\frac{k(k+3)}{2}}}{(1 - x^3)(1 - x^4) \dots (1 - x^k)}. \tag{2.9}$$

In terms of integer partitions, we obtain the following proposition.

**Proposition 2.1:** *For any integer $n \geq 9$, the number of partitions with odd parts larger or equal to 5 is equal to the number of partitions with distinct parts, at least three parts larger or equal to 2, up to two 1 part (the part 1 may not occur, or occur once, or occur twice), and the three largest parts consecutive.*

$\square$

(Proposition 2.1, here, is also <span style="color:red">Proposition 2.1</span> in <span style="color:red">[H]</span>, p. 4.)

The cyclotomic polynomial $1 + x + x^2$ in (2.9) obviously signifies that, adding up three consecutive terms of the butterfly sequence (see (1.10)), gives the number of partitions with odd parts larger or equal to 5, as we can see for the first few terms by direct computation:

$$\{t(n)\} = \{s(n) + s(n-1) + s(n-2)\} = \{\, 1, 0, 0, 0, 0, 1, 0, 1, 0, 1, 1, 1, 1, 1, 2, \tag{2.10}$$

$$2, 2, 3, 3, 4, 4, 5, 5, 6, 7, 8, 9, 10, 11, 13, , 14, 16, 18, 20, 23, 25, 28, 31, 35, 39,$$

$$41, 46, 51, 59, 66, 72, 80, 88, 98, 108, 119, \dots \,\},$$

for $n \geq 0$ (we will get back to this sequence at the end of the present paper).

Furthermore, if we note that the polynomial with degree smaller or equal than 8 in (2.9) is obviously $1 + x^5 + x^7$, we can write the generating function identity for the partitions of n with odd parts larger or equal to 5 as

$$\prod_{h=2}^{\infty} \frac{1}{1 - x^{2h+1}} = 1 + x^5 + x^7 + \sum_{k=2}^{\infty} \frac{(1 + x + x^2) x^{\frac{k(k+3)}{2}}}{(1 - x^3)(1 - x^4) \dots (1 - x^k)}. \tag{2.11}$$

For the complete generating function of the butterfly sequence, finally, we just divide both sides of (2.11) by $1 + x + x^2$ and, as

$$\frac{1 + x^5 + x^7}{1 + x + x^2} = 1 - x + x^3 - x^4 + x^5,$$

we see that (2.11) results in a complete filtered (all terms, for any nonnegative integer n, infinite product and filtration by number of parts series) generating function identity for the butterfly sequence:

$$\frac{1}{1 + x + x^2} \prod_{h=2}^{\infty} \frac{1}{1 - x^{2h+1}} = 1 - x + x^3 - x^4 + x^5 +$$

$$+ \sum_{k=2}^{\infty} \frac{x^{\frac{k(k+3)}{2}}}{(1 - x^3)(1 - x^4) \dots (1 - x^k)}. \tag{2.12}$$

The identities (2.11) and (2.12) also clarify the exceptions in the second difference pentagonal number structure of the sequence (1.10) (see Remark 1.2, here).

## 3 THE BUTTERFLY SEQUENCE AS A SEQUENCE OF PARTITIONS WITH ODD PARTS LARGER OR EQUAL TO 3

There are several procedures to pair partitions with distinct parts with partitions with odd parts. The one we chose to apply in this section to the butterfly partitions is based on Euler's versatile "splitting and merging" procedure, and on the two special features of the butterfly partition. The "splitting and merging" procedure, here, will consist of splitting most of the even parts in the partitions with distinct parts to produce partitions with odd parts, and, vice versa, merging some of the odd parts into even parts (in the partitions with odd parts to produce partitions with distinct parts). The first special feature of the butterfly partition, the presence of the three largest consecutive parts, will be of use, here, in producing three similar (odd) parts in the partitions with odd parts corresponding to the butterfly partitions. The other special feature, "the absence of 1 parts" in the butterfly partitions, will be interpreted as "the sum of the parts that are nonnegative integral powers of two is an even number". This even number, in turn, will be added to the largest part of the partition with odd parts.

To better understand the point of "adding the nonnegative integral powers of two to the largest part", the reader may find useful to recall some of the details of the procedure that leads to the interpretation of the sequence (1.2) as the sequence of the partitions of $n$, for $n \geq 4$, with distinct parts that are not powers of two. This interpretation can be obtained, partition by partition, by adding the sum of the (distinct) powers of two to one of the two largest parts (the largest, or, if the result of the sum is itself a power of two, the second largest). It may also be useful to recall the generating function interpretation of (1.2), the fact that the left side of (2.3) can be written as

$$(1-x)\prod_{n=1}^{\infty}(1+x^n) = (1-x^2)\prod_{n=2}^{\infty}(1+x^n) = (1-x^4)\prod_{n=3}^{\infty}(1+x^n)$$

$$= (1-x^8)\prod_{n=3,n\neq4}^{\infty}(1+x^n) = \cdots = (1-x^{2m})\prod_{n=1,n\neq2^k,k=1,2,\dots,m-1}^{\infty}(1+x^n),$$

and, as $x$ is a formal variable (or, as usual, for $x$ a complex variable in the interval of convergence $|x| < 1$),

$$(1-x)\prod_{n=1}^{\infty}(1+x^n) = \prod_{n=1,n\neq2^k,k=1,2,3,\dots}^{\infty}(1+x^n).$$

With these special features in mind, we begin our visualization of the butterfly sequence as a sequence of partitions with odd parts greater than or equal to 3. To focus on the even and odd parts among the three largest parts of the partitions, we will use the following definition.

**Definition 3.1:** For any positive integer $n$, $n \geq 6$, we define two subsequences of the butterfly sequence. The first, the sequence $\{s_e(n)\}$ as the number of butterfly partitions of $n$ with *the second largest part even,* that is the number of partitions of the form

$$p_1 = 2m+1 > p_2 = 2m > p_3 = 2m-1 > p_4 > \cdots > p_k \geq 2, \tag{3.1}$$

with

$$\sum_{i=1}^{k} p_i = n,$$

and $k \geq 3$. And the second subsequence, also for $n \geq 6$, the sequence $\{s_o(n)\}$ of the number of butterfly partitions of $n$ with *the second largest part odd,* that is the number of partitions of the form

$$p_1 = 2m > p_2 = 2m-1 > p_3 = 2m-2 > p_4 > \cdots > p_k \geq 2, \tag{3.2}$$

with

$$\sum_{i=1}^{k} p_i = n,$$

and $k \geq 3$.

□

(Definition 3.1, here, is also Definition 3.1 in [H], p. 5.)

The first few terms of the sequences $\{s_e(n)\}$ and $\{s_o(n)\}$ can be quickly computed as in

$$\{s_e(n)\} = \{0, 0, 0, 0, 0, 0, 1, 0, 1, 0, 0, 0, 1, 0, 1, 1, 1, 1, 2, 1, 1, 2, 1, 2, 3, 2, 3, 4, 3, 5, 5, 5, 6,$$
$$7, 6, 7, 9, 9, 11, 13, 12, 15, 17, 17, 20, 23, \dots \}, \qquad (3.3)$$

and, respectively, as in

$$\{s_o(n)\} = \{0, 0, 0, 1, 0, 0, 0, 0, 0, 1, 0, 1, 1, 0, 1, 1, 0, 1, 1, 2, 2, 2, 2, 3, 2, 3, 4, 3, 4, 5, 4, 6,$$
$$7, 7, 7, 10, 9, 11, 13, 12, 15, 17, 17, 20, 22, \dots \}, \qquad (3.4)$$

both (3.3) and (3.4) for $n = 6, 7, 8 \dots$.

We will refer to the sequence $\{s_e(n)\}$, for $n \geq 6$, as the *butterfly sequence with the second largest part even*, and to the sequence $\{s_o(n)\}$, for $n \geq 6$, as the *butterfly sequence with the second largest part odd*.

**Remark 3.2:** We choose to start the sequences (3.3) and (3.4) at $n = 6$, although the smallest nonzero such partition with the second part even occurs for $n = 12$ ($p_1 = 5 > p_2 = 4 > p_3 = 3$), while the smallest nonzero such partition with the second part odd occurs for $n = 9$ ($p_1 = 4 > p_2 = 3 > p_3 = 2$), so that, obviously, in both (3.1) and (3.2) $m \geq 2$.

□

(For the readers of [H], note that Remark 3.2, here, does not correspond to Remark 3.2 in [H]. Remark 3.3, below, is an extended version of Remark 3.2 in [H].)

For both kinds (second part even or odd) of butterfly partitions in the above definition, we also denote by $t$ half the sum of all the parts *among $p_4, p_5, \dots p_k$* in (3.1), or, respectively, in (3.2), that are distinct, positive, and integral powers of two

$$t = \tfrac{1}{2} \sum_{j=1}^{h} p_{i_j} \qquad (3.5)$$

denoting $p_{i_j} = 2^{i_j}, j = 1, 2, \dots, h,$ all the parts (in decreasing order, as usual) in (3.1) or (3.2) among $p_4, p_5, \dots p_k$ that are powers of two, so that

$$2t = \sum_{j=1}^{h} p_{i_j}$$

is the sum of all the parts that are powers of two.

The construction of partitions with odd parts greater than or equal to 3 corresponding to the partitions (3.1) and, respectively, (3.2) is presented in this section in four parts.

**Step I.** For each butterfly partition in (3.1) (with the second part, $p_2 = 2m$, even), we construct a specific partition of the same integer n, for $n \geq 6$ (see Remark 3.1, above),

$$q_1 \geq q_2 = q_3 \geq q_4 \geq \ \dots \ \geq q_r = 3 \tag{3.6}$$

with

$$\sum_{j=1}^{r} q_j = n,$$

and each part $q_j, j = 1, 2, \dots, r$ odd, greater than or equal to 3, $r \geq 4$ (the partition (3.6) must have at least four parts), the second largest part equal to the third, and *the third largest part greater than the fourth largest unless the second, the third, and the fourth largest parts are all equal to 3*, as following.

We set the three largest parts

$$q_1 = 2m - 1 + 2t, \ q_2 = 2m - 1, \ q_3 = 2m - 1, \tag{3.7}$$

and we set the smallest part

$$q_r = 3 \tag{3.8}$$

(note that $q_1 + q_2 + q_3 + q_r$ is equal to the sum of $p_1 + p_2 + p_3$, as in (3.1), and $2t$, as in (3.5)). As for all the remaining parts $q_4 \geq q_5 \geq \ \dots \ \geq q_{r-1}$, they are formed with all the odd parts among $p_4 > p_5 > \dots > p_k$ in (3.1) together with all the odd parts obtained by the usual Euler's splitting procedure applied to all the even parts among $p_4 > \dots > p_k$ in (3.1) that are not powers of two. Note that, in this scenario (the scenario with the second part, $p_2 = 2m$ in (3.1), even), the corresponding second and the third largest (odd) parts in (3.6) are equal, and note that the fourth largest part in (3.6) is necessarily smaller than the third largest part unless the second, the third, and the fourth largest parts are all equal to 3.

**Example 1:** The butterfly partition

$$p_1 = 5 > p_2 = 4 > p_3 = 3,$$

(a partition with $m = 2$ and $k = 3$ in (3.1)) corresponds to the partition with four parts, all equal to 3, as in (3.6)

$$q_1 = q_2 = q_3 = q_4 = 3,$$

($m = 2$, $t = 0$, and $r = 4$ in (3.6), (3.7), and (3.8)).

□

**Example 2:** The butterfly partition of $n = 44$

$$p_1 = 9 > p_2 = 8 > p_3 = 7 > p_4 = 6 > p_5 = 5 > p_6 = 4 > p_8 = 3 > p_8 = 2$$

(a partition with $m = 4$ and $k = 8$ in (3.1)) corresponds to the partition with odd parts greater than or equal to 3

$$q_1 = 13 > q_2 = q_3 = 7 > q_4 = 5 > q_5 = q_6 = q_7 = q_8 = 3$$

($m = 4$, $t = 2$, and $r = 8$ in (3.6), (3.7), and (3.8)). Note that $q_1 = 7 + 2t = 7 + p_6 + p_8$, $q_7 = p_7$, and $q_5 = q_6 = 3$ are obtained by splitting $p_4 = 6$.

□

**Step II.** For each butterfly partition in (3.2) (with the second part, $p_2 = 2m - 1$, odd), we construct a specific partition of the same integer n, for $n \geq 6$,

$$q_1 > q_2 > q_3 = q_2 - 2 \geq q_4 \geq \ \dots \ \geq q_r \geq 3 \qquad (3.9)$$

with

$$\sum_{j=1}^{r} q_j = n,$$

and each part $q_j$, $j = 1, 2, \dots, r$ odd, greater than or equal to 3, $r \geq 3$ (the partition (3.9) must have at least three parts), the largest part *at least* two units greater than the second largest part, and the second largest part *exactly* two units greater than to the third, as following.

For $m = 2$, we associate the butterfly partition

$$p_1 = 2m = 4 > p_2 = 2m - 1 = 3 > p_3 = 2m - 2 = 2,$$

to the partition with odd parts

$$q_1 = q_2 = q_3 = 2m - 1 = 3$$

(which is the *initial partition* with odd parts of Step II). For every other butterfly partition in (3.2), set

$$q_1 = 2m + 1 + 2t, \ q_2 = 2m - 1, \ q_3 = 2m - 3 \qquad (3.10)$$

(note that $q_1 + q_2 + q_3$ is equal to the sum of $p_1 + p_2 + p_3$, as in (3.2), and $2t$, as in (3.5)). As for all the remaining parts $q_4 \geq q_5 \geq \ ... \ \geq q_r$, as in Step I, they are formed with all the odd parts among $p_4 > p_5 > \cdots > p_k$ in (3.2) together with all the odd parts obtained by the usual Euler's splitting procedure applied to all the even parts among $p_4 > \cdots > p_k$ in (3.2) that are not powers of two. Note that, in this scenario (the scenario with the second part, $p_2 = 2m - 1$ in (3.2), odd), *except for the initial partition,* the largest part in (3.9) is set to be at least two units greater than the second largest, and the second largest part in (3.9) is set to be two units greater than the third largest part.

**Example 3:** The butterfly partition of $n = 52$

$$p_1 = 10 > p_2 = 9 > p_3 = 8 > p_4 = 7 > p_5 = 6 > p_6 = 5 > p_7 = 4 > p_8 = 3$$

($m = 4$ and $k = 8$ in (3.2)) corresponds to the partition with odd parts greater than or equal to 3

$$q_1 = 15 > q_2 = 9 > q_3 = q_4 = 7 > q_5 = 5 > q_6 = q_7 = q_8 = 3$$

($m = 4$, $t = 2$, and $r = 8$ in (3.6), (3.7), and (3.8)). Note that $q_1 = 11 + p_7$, $q_2 = p_2 = 9$, $q_3 = p_3 - 1 = q_4 = p_4 = 7$, $q_5 = p_6$, $q_6 = q_7 = 3$ are obtained by splitting $p_4 = 6$, $q_7 = p_8 = 3$, and note that the part $p_3 = 8$, although a power of two, is not incorporated into the first part of (3.9) because is also one of the three largest parts of (3.2).

□

**Step III.** For each partition with odd parts greater than or equal to 3 of the form (3.6), with the three largest parts as in (3.7), and the smallest part as in (3.9), or, respectively, of the form (3.9), with (3.10), and, hence, with at least three parts, under specific *merging conditions*, we reconstruct the butterfly partitions (3.1), and, respectively, (3.2).

The procedures of Steps I and II, above, are obviously reversible. The parts $q_1, q_2, q_3$, and $q_r$ in (3.6), or respectively, $q_1, q_2,$ and $q_3$ in (3.9) are easily deconstructed into the parts $p_1, p_2, p_3$, and all the parts in (3.1), or, respectively, (3.2) that are distinct nonzero powers of two (as the number $2t$ in (3.5) is identified in a sum of distinct nonzero powers of two). The remaining odd parts $q_4, q_5,$ and $q_{r-1}$ in (3.6), or, respectively, $q_4, q_5,$ and $q_r$ in (3.9) are merged into distinct parts in (3.1) and (3.2), as, for each odd number $q$ greater than or equal to 3, the number of the remaining odd parts in (3.6) and (3.9) equal to $q$ is written as a sum of distinct powers of two, and each of these powers of two, multiplied by $q$, produces a (distinct) part of (3.1), or, respectively, (3.2).

Since each partition in (3.6) is distinct from each partition in (3.9) (for the second and third largest parts equal versus the second largest part two units larger than the third), the above procedures give a one-to-one correspondence between the butterfly partitions and some of the partitions with odd parts greater than or equal to 3. To see this through, we need to compute precisely the *merging cap conditions* that determine the feasible sizes of both $2t$ in (3.5) and the number of repeated odd parts of the partitions in (3.6) and (3.9).

We begin with the smallest partition with at least three parts, each part odd and greater than or equal to 3, $q_1 = q_2 = q_3 = 3$, and, for this partition, we construct the butterfly partition $p_1 = 4 > p_2 = 3 > p_3 = 2$. Then, consider a partition of the form (3.6), or, respectively, (3.9). To see what we call the *merging cap information of the partition,* we begin by reconstructing the parts in (3.1), and, respectively, (3.2) that are powers of two, from the difference between the two largest parts of the partition. If the difference between the two largest parts in (3.6) is zero, or, respectively, if the difference between the two largest parts in (3.9) is two, in other words, if $t = 0$ (so that the corresponding butterfly partition will have no parts that are powers of two), we say that the partition with odd parts has no merging cap on the powers of two (as the corresponding butterfly partition has no parts, among $p_4 > p_5 \dots > p_k$ that are powers of two). For all the other partitions of the form (3.6), or, respectively (3.9), $t > 0$, and we can write, first, for the partitions of the form (3.6), the difference between the two largest parts as

$$q_1 - q_2 = 2t = 2 \sum_{i=0}^{k} 2^{r_{2,i}}, \tag{3.11}$$

($t$ in binary form), and, second, using the same notation, for the partitions of the form (3.9), the similar difference between the largest part $q_1$ and the number $q_2 + 2$ as

$$q_1 - 2 - q_2 = 2t = 2 \sum_{i=0}^{k} 2^{r_{2,i}}, \tag{3.12}$$

so that, in both (3.11) and (3.12), the exponents are in decreasing order

$$r_{2,0} > r_{2,1} > \cdots > r_{2,k} \geq 0,$$

and, in particular, $2^{r_{2,0}}$ denotes the largest power of 2 in the binary form of the integer $t > 0$. Then, the merging cap embodied in the two largest part of (3.6), or, respectively, (3.9) is the condition that the largest nonzero power of 2 in $t$ in (3.11) recreates a part $2 \cdot 2^{r_{2,0}}$ in (3.2) no larger than $p_3 - 1$ (the largest possible value of $p_4$ in a butterfly partition) and, therefore, requires

$$2 \cdot 2^{r_{2,0}} \leq q_3 - 1, \text{ or, equivalently, } r_{2,k} \leq (\log_2 (q_3 - 1)) - 1, \qquad (3.13)$$

for $q_3 (= 2m - 1 = p_3)$ as in (3.6) and (3.7), and, respectively,

$$2 \cdot 2^{r_{2,0}} \leq q_3, \text{ or, equivalently, } r_{2,k} \leq (\log_2 q_3) - 1, \qquad (3.14)$$

for $q_3 (= 2m - 3 = p_3 - 1)$ as in (3.9) and (3.10).

Similar merging caps apply to the number of repeated odd parts among the parts

$$q_4 \geq q_5 \geq \ldots \geq q_r \qquad (3.15)$$

in (3.6) and (3.9), as follows. For any odd integer $q > 3$, if $q \neq q_i, i = 4, 5, \ldots, r$, there is no merging cap (as the corresponding butterfly partition will have no parts that are the product of $q$ and a nonnegative integral power of two). If, on the other hand $q = q_i$, for some $i = 4, 5, \ldots, r$, we write the number of parts equal to $q$ (among the parts in (3.15)) in binary form as

$$u = \sum_{i=0}^{k} 2^{r_{q,i}} \qquad (3.16)$$

so that, in (3.16), the exponents are in decreasing order

$$r_{q,0} > r_{q,1} > \cdots > r_{q,k} \geq 0,$$

and, in particular, $2^{r_{q,0}}$ denotes the largest power of 2 in the binary form of the integer $u > 0$. Then, the merging cap on $q$ is the condition that the largest nonzero power of 2 in $u$ in (3.16) recreates (through the merging) a part $q \cdot 2^{r_{q,0}}$ in (3.2) no larger than $p_3 - 1$ (the largest possible value of $p_4$ in a butterfly partition) and, therefore, requires

$$q \cdot 2^{r_{q,0}} \leq q_3 - 1, \text{ or, equivalently, } r_{q,k} \leq (\log_2 (q_3 - 1)) - (\log_2 q) \qquad (3.17)$$

for $q_3 (= 2m - 1 = p_3)$ as in (3.6) and (3.7), and, respectively,

$$q \cdot 2^{r_{q,0}} \leq q_3, \text{ or, equivalently, } r_{q,k} \leq (\log_2 q_3) - (\log_2 q) \qquad (3.18)$$

for $q_3 (= 2m - 3 = p_3 - 1)$ as in (3.9) and (3.10).

For $q = 3$, there are no merging cap, if there is only one part, $q_r$ in (3.15), equal to 3 in the (3.6) scenario, and, respectively, if there are no parts $q_i, i = 4, 5, \ldots, r$ in (3.15) equal to 3 in the (3.9) scenario.

If, on the other hand, there is more than one part $q_i, i = 4, 5, \ldots, r$ in (3.15) equal to 3 in the (3.6) scenario, or, respectively, if there is at least one part $q_i, i = 4, 5, \ldots, r$ in (3.15) equal to 3 in the (3.9) scenario, we write both the number of parts *minus one* equal to 3 in (3.15) *in the* (3.6) *scenario*, and the number of parts equal to 3 in (3.15) in the (3.9) scenario, in the same binary form as

$$v = \sum_{i=0}^{k} 2^{r_{3,i}} \tag{3.19}$$

so that, in (3.16), the exponents are in decreasing order

$$r_{3,0} > r_{3,1} > \cdots > r_{3,k} \geq 0,$$

and, in particular, $2^{r_{3,0}}$ denotes the largest power of 2 in the binary form of the integer $v > 0$. Then, the merging cap on the number of parts equal to 3 is the condition that the largest nonzero power of 2 in $v$ in (3.19) recreates (through the merging) a part $3 \cdot 2^{r_{3,0}}$ in (3.2) no larger than $p_3 - 1$ (the largest possible value of $p_4$ in any partition with distinct parts of the form (3.2)) and, therefore, requires

$$3 \cdot 2^{r_{3,0}} \leq q_3 - 1, \text{ or, equivalently, } r_{3,k} \leq \left(\log_2 (q_3 - 1)\right) - (\log_2 3) \tag{3.20}$$

for $q_3 (= 2m - 1 = p_3)$ as in (3.6) and (3.7), and, respectively,

$$3 \cdot 2^{r_{3,0}} \leq q_3, \text{ or, equivalently, } r_{3,k} \leq \left(\log_2 q_3\right) - (\log_2 3) \tag{3.21}$$

for $q_3 (= 2m - 3 = p_3 - 1)$ as in (3.9) and (3.10).

**Example 4:** The partition of $n = 14$ with odd parts

$$q_1 = 5 > q_2 = q_3 = q_4 = 3$$

is the second smallest (in the lexicographic order) partition as in (3.6), and is easily deconstructed as

$$q_1 = 2m - 1 + 2t > q_2 = q_3 = q_4 = 2m - 1$$

with $m = 2$, and $t = 1$, and, therefore associated to the butterfly partition

$$p_1 = 2m + 1 = 5 > p_2 = 2m = 4 > p_3 = 2m - 1 = 3 > p_4 = t = 2.$$

Note that the only merging cap applicable, (3.13), in this example is barely satisfied (the partition reaches the merging cap), as the corresponding butterfly partition is semi-triangular (triangular except for the part 1), or, equivalently, in terms of (3.13)

$r_{2,0} = 0 \leq \left(\log_2 (3 - 1)\right) - 1 = 0.$

Note that the merging cap (3.13) tells us that the only partitions of the form (3.6) with four parts, and $q_2 = q_3 = q_4 = 3$ are the partitions with odd parts in examples 1 and 4 ($q_1 = 3$ or $q_1 = 5$).

□

**Example 5:** The partition of $n = 90$ with twelve odd parts

$$q_1 = 25 > q_2 = q_3 = 11 > q_4 = 9 > q_5 = 7 > q_6 = q_7 = q_8 = 5 >$$

$$> q_9 = q_{10} = q_{11} = q_{12} = 3$$

is deconstructed first observing that $m = 6$ (for $q_2 = q_3 = 2m - 1 = 11$, and $q_{12} = 3$) and $t = 7$ (for $q_1 = 2m - 1 + 2t = 25$). This gives the first three parts of the corresponding butterfly partition with the second part even

$$p_1 = 2m + 1 = 13 > p_2 = 2m = 12 > p_3 = 2m - 1 = 11$$

and, since $2t = 14 = 8 + 4 + 2$, all the parts of the corresponding butterfly partition that are distinct powers of two

$$p_{i_1} = 8 > p_{i_2} = 4 > p_{i_3} = 2$$

(the parts $p_6, p_{10}$, and $p_{12}$, below). Second, observe that single odd parts $q_4 = 9$ and $q_5 = 7$ will translate intact into the resulting butterfly partition (the parts $p_5$ and $p_7$), while the multiple odd parts $q_6 = q_7 = q_8 = 5$, and, respectively $q_9 = q_{10} = q_{11} = 3$, will merge into $p_4 = 2 \cdot 5 = 10$, and $p_9 = 5$ (for there are $2 + 1$ parts equal to 5 in the given partition with odd parts), and respectively, into $p_8 = 2 \cdot 3 = 6$, and $p_{11} = 3$ (for, again, there are $2 + 1$ parts equal to 3 in the partition with odd parts, not including the part $q_{12} = 3$ which is already incorporated in the first three parts of the result, the parts $p_1 = 13 > p_2 = 12 > p_3 = 11$ of the butterfly partition below). The result is the following butterfly partition of $n = 90$

$$p_1 = 13 > p_2 = 12 > p_3 = 11 > p_4 = 10 > p_5 = 9 > p_6 = 8 > p_7 = 7 > p_8 = 6 >$$
$$> p_9 = 5 > p_{10} = 4 > p_{11} = 3 > p_{12} = 2.$$

As in the previous example, all the merging caps applicable, (3.13), (3.17) with $q = 5$, and (3.20) are barely satisfied (the resulting butterfly partition is semi-triangular):

$r_{2,0} = 3 \leq (\log_2 (11 - 1)) - 1$ ( $r_{2,0} = 3$ is the integral part of $(\log_2 10) - 1$,

$r_{5,0} = 1 \leq (\log_2 (11 - 1)) - (\log_2 5) = 1$,

and

$r_{6,0} = 1 \leq (\log_2 (6 - 1)) - (\log_2 3) = 1$.

which tells us that the only partitions of the form (3.6) with four parts, and $q_2 = q_3 = q_4 = 3$ are the partitions with odd parts in examples 1 and 4 ($q_1 = 3$ or $q_1 = 5$).

$\square$

**Remark 3.3:** Note that the procedures in Step I and Step II can be switched in the following way.

The procedure in Step I that applies to (3.1), and produces (3.6), (3.7) and (3.8), can be replaced with a procedure, similar to the procedure in Step II, that applies to (3.1), and that uses the same notation introduced in this section. For $m = 2$ in (3.1), for the butterfly partition

$$p_1 = 2m + 1 = 5 > p_2 = 2m = 4 > p_3 = 2m - 1 = 3,$$

we construct the partition

$$q_1 = q_2 = q_3 = q_4 = 2m - 1 = 3$$

(the *initial partition* with odd parts of Step I). As for all the remaining partitions in (3.1), in place of (3.6), (3.7) and (3.8), the corresponding partitions with odd parts are

$$q_1 > q_2 > q_3 = q_2 - 2 \geq q_4 \geq \ ... \ \geq q_r = 3 \tag{3.22}$$

with the three largest parts

$$q_1 = 2m + 1 + 2t, \ q_2 = 2m - 1, \ q_3 = 2m - 3, \tag{3.23}$$

and the smallest part

$$q_r = 3 \tag{3.24}$$

incorporating the three largest parts of (3.1) together with the parts among $p_4 > p_5 > \cdots > p_k$ that are powers of two, and the same splitting procedure of Step I and Step II applied to the even parts among the remaining parts $p_4 > p_5 > \cdots > p_k$ in (3.1), to produce the parts $q_4 \geq q_5 \geq \ ... \ \geq q_{r-1}$ (compare (3.6), (3.7) and (3.8) with (3.22), (3.23) and (3.24)).

On the other hand, to keep the partitions with odd parts corresponding to the partitions in (3.1) distinct from the partitions with odd parts corresponding to (3.2), if we change the procedure that applies to (3.1), we must change the procedure that applies to (3.2). For this, the procedure in Step II that applies to the partitions in (3.2) and produces (3.9) and (3.10), is replaced with a procedure, similar to Step I, that applies to (3.2), using the same notation, and produces in place of (3.9) and (3.10),

$$q_1 \geq q_2 = q_3 \geq q_4 \geq \ ... \ \geq q_r = 3 \tag{3.25}$$

with the parts

$$q_1 = 2m - 1 + 2t, \ q_2 = 2m - 1, \ q_3 = 2m - 1, \tag{3.26}$$

incorporating the three largest parts of (3.2) together with the parts among $p_4 > p_5 > \cdots > p_k$ that are powers of two, and the same splitting procedure of Step I and Step II applied to the even parts among the remaining parts $p_4 > p_5 > \cdots > p_k$ in (3.1), to produce the parts $q_4 \geq q_5 \geq \ldots \geq q_r$ (compare (3.9) and (3.10) with (3.25) and (3.26)).

If we switch the procedures in Step I and Step II, the merging caps (3.13), (3.14), (3.17), (3.18), (3.20), and, respectively, (3.21), *if applicable*, need to be adjusted, and, using the same notation, replaced with (3.27), (3.28), (3.29), (3.30), (3.31), and, respectively, (3.32).

The merging cap embodied in the two largest part of (3.22) requires

$$2 \cdot 2^{r_{2,0}} \leq q_3 + 1, \text{ or, equivalently, } r_{2,k} \leq (\log_2 (q_3 + 1)) - 1, \qquad (3.27)$$

for $q_3 (= 2m - 3 = p_3 - 2)$ as in (3.22) and (3.23), and, respectively, the merging cap for the two largest parts in (3.25)

$$2 \cdot 2^{r_{2,0}} \leq q_3 - 2, \text{ or, equivalently, } r_{2,k} \leq (\log_2 (q_3 - 2)) - 1, \qquad (3.28)$$

for $q_3 (= 2m - 1 = p_3 + 1)$ as in (3.25) and (3.26).

The merging cap on the number of parts equal to any odd number $q \neq 3$ among the parts $q_4 \geq q_5 \geq \ldots \geq q_{r-1}$ in (3.22), requires

$$q \cdot 2^{r_{q,0}} \leq q_3 + 1, \text{ or, equivalently, } r_{q,k} \leq (\log_2 (q_3 + 1)) - (\log_2 q) \qquad (3.29)$$

for $q_3 (= 2m - 3 = p_3 - 2)$ as in (3.22) and (3.23), and, respectively, the merging cap on the number of parts equal to any odd number $q \neq 3$ among the parts $q_4 \geq q_5 \geq \ldots \geq q_r$ in (3.25), requires

$$q \cdot 2^{r_{q,0}} \leq q_3 - 2, \text{ or, equivalently, } r_{q,k} \leq (\log_2 (q_3 - 2)) - (\log_2 q) \qquad (3.30)$$

for $q_3 (= 2m - 1 = p_3 + 1)$ as in (3.25) and (3.26).

Finally, for the merging cap on the number of parts equal to 3, recall the notation (3.19) which requires to set one part equal to 3 aside in the partitions (3.22), (3.23) and (3.24). For this, the merging cap on the number of parts equal to 3 among the parts $q_4 \geq q_5 \geq \ldots \geq q_{r-1}$ in (3.22), requires

$$3 \cdot 2^{r_{3,0}} \leq q_3 + 1, \text{ or, equivalently, } r_{3,k} \leq (\log_2 (q_3 + 1)) - (\log_2 3) \qquad (3.31)$$

for $q_3 (= 2m - 3 = p_3 - 2)$ as in (3.22) and (3.23), and, respectively, the merging cap on the number of parts equal to 3 among the parts $q_4 \geq q_5 \geq \ldots \geq q_r$ in (3.25),

$$3 \cdot 2^{r_{3,0}} \leq q_3 - 2, \text{ or, equivalently, } r_{3,k} \leq (\log_2 (q_3 - 2)) - (\log_2 3) \qquad (3.32)$$

for $q_3 (= 2m - 1 = p_3 + 1)$ as in (3.25) and (3.26).

□

(For the readers of [H], Remark 3.3, here, is an extended version of Remark 3.2 in [H], p. 7.)

**Example 6:** The butterfly partition of $n = 27$

$$p_1 = 7 > p_2 = 6 > p_3 = 5 > p_4 = 4 > p_5 = 3 > p_6 = 2$$

(a partition with $m = 3$ and $k = 6$ in (3.1)), *using the procedure in* Step I, above, corresponds to the partition with odd parts greater than or equal to 3

$$q_1 = 11 > q_2 = q_3 = 5 > q_4 = q_5 = 3$$

($m = 3$, $t = 6$, and $r = 5$ in (3.6), (3.7), and (3.8)), while *using the switched procedure in* Remark 3.3, here, corresponds to the partition with odd parts greater than or equal to 3

$$q_1 = 13 > q_2 = 5 > q_3 = q_4 = q_5 = 3$$

($m = 3$, $t = 6$, and $r = 5$ in (3.22), (3.23), and (3.24)).

□

**Example 7:** The butterfly partition of $n = 27$

$$p_1 = 10 > p_2 = 9 > p_3 = 8$$

(a partition with $m = 5$ and $k = 3$ in (3.2)), *using the procedure in* Step II, above, corresponds to the partition with odd parts greater than or equal to 3

$$q_1 = 11 > q_2 = 9 > q_3 = 7$$

($m = 5$, $t = 0$, and $r = 3$ in (3.9) and (3.10)), while *using the switched procedure in* Remark 3.3, here, corresponds to the partition with odd parts greater than or equal to 3

$$q_1 = 9 > q_2 = 9 > q_3 = 9$$

($m = 3$, $t = 0$, and $r = 3$ in (3.25), and (3.26)).

□

The following proposition summarizes the results of Steps I, II, and III in this section.

**Proposition 3.4:** *(a) For any positive integer $n$, $n \geq 6$, let $\{o_e(n)\}$ be the number of partitions of the form*

$$q_1 \geq q_2 = q_3 \geq q_4 \geq \ ... \ \geq q_r = 3$$

*with*

$$\sum_{j=1}^{r} q_j = n,$$

and each part $q_j, j = 1, 2, \ldots, $ r odd, greater than or equal to 3, r ≥ 4 (the partition must have at least four parts), the second largest part equal to the third, and the third largest part greater than the fourth largest unless the second, the third, and the fourth largest parts are all equal to 3, subject to the following conditions (merging caps):

*(i)*      $2 \cdot 2^{r_{2,0}} \leq q_3 - 1,$

*(ii)*      $q \cdot 2^{r_{q,0}} \leq q_3 - 1,$ and

*(iii)*      $3 \cdot 2^{r_{3,0}} \leq q_3 - 1,$

where $2^{r_{2,0}}, 2^{r_{q,0}},$ and, respectively, $2^{r_{3,0}},$ denote the largest powers of two less than or equal to the three quantities $2t, u,$ and, respectively, $v,$ with $2t = q_1 - q_2$ (the difference between the two largest parts of the partition), $u = u(q)$ the number of parts, among $q_4, q_5, \ldots q_r,$ that are equal to an odd integer $q,$ for any odd integer $q \neq 3,$ and, respectively, $v$ the number of parts, among $q_4, q_5, \ldots q_{r-1},$ that are equal to 3. (If any of these quantities is zero, the corresponding condition (i), (ii), or (iii) does not apply.) Then, for any $n \geq 6,$

$o_e(n) = s_e(n)$ *(see Definition 3.1).*

*(b) For any positive integer $n, n \geq 6,$ let $\{o_o(n)\}$ be the number of partitions of the form*

$$q_1 = q_2 = q_3 = 3$$

*for* n $= q_1 + q_2 + q_3 = 9,$ *or of the form*

$$q_1 > q_2 > q_3 = q_2 - 2 \geq q_4 \geq \ldots \geq q_r \geq 3$$

*with*

$$\sum_{j=1}^{r} q_j = n,$$

and each part $q_j, j = 1, 2, \ldots, $ r odd, greater than or equal to 3, r ≥ 3 (the partition must have at least three parts), the largest part at least two units greater than the second largest part,

$$q_1 - q_2 \geq 2,$$

and the second largest part exactly two units greater than the third largest part, subject to the following conditions (merging caps):

*(i)*      $2 \cdot 2^{r_{2,0}} \leq q_3,$

*(ii)*      $q \cdot 2^{r_{q,0}} \leq q_3,$ and

*(iii)*     $3 \cdot 2^{r_{3,0}} \le q_3,$

*where $2^{r_{2,0}}$, $2^{r_{q,0}}$, and, respectively, $2^{r_{3,0}}$, denote the largest powers of two less than or equal to the three quantities $2t$, $u$, and, respectively, $v$, with $2t = q_1 - q_2 - 2$ (the difference between the largest part of the partition and the second largest plus two units), $u = u(q)$ the number of parts, among $q_4$, $q_5$, ... $q_r$, that are equal to an odd integer $q$, for any odd integer $q \ne 3$, and, respectively, $v$ the number of parts, among $q_4$, $q_5$, ... $q_{r-1}$, that are equal to 3. (If any of these quantities is zero, the corresponding condition (i), (ii), or (iii) does not apply.)*

*Then, for any $n \ge 6$,*

$\mathrm{o_o}(n) = \mathrm{s_o}(n)$ *(see Definition 3.1).*

*(c) For any positive integer $n$, $n \ge 6$,*

$$\mathrm{o_e}(n) + \mathrm{o_o}(n) = \mathrm{s}(n)$$

*(using the same notation of parts (a) and (b), and Definition 1.1).*

□

(Proposition 3.4, here, corresponds to <span style="color:red">Proposition 3.3</span> in <span style="color:red">[H]</span>, pp. 7 - 8.)

Using Remark 3.3, we obtain the following slight variation of Proposition 3.4.

**Proposition 3.5:** *(a) For any positive integer $n$, $n \ge 6$, let $\{\mathrm{o'_e}(n)\}$ be the number of partitions of the form*

$$q_1 = q_2 = q_3 = q_4 = 3$$

*for $\mathrm{n} = q_1 + q_2 + q_3 + q_4 = 12$, or of the form*

$$q_1 > q_2 > q_3 = q_2 - 2 \ge q_4 \ge \ ... \ge q_r \ge 3$$

*with*

$$\sum_{j=1}^{r} q_j = n,$$

*and each part $q_j$, $j = 1, 2, ... , \mathrm{r}$ odd, greater than or equal to 3, $\mathrm{r} \ge 4$ (the partition must have at least four parts), the largest part at least two units greater than the second largest part,*

$$q_1 - q_2 \ge 2,$$

*and the second largest part exactly two units greater than the third largest part, subject to the following conditions (merging caps):*

*(i)*      $2 \cdot 2^{r_{2,o}} \leq q_3 + 1,$

*(ii)*     $q \cdot 2^{r_{q,o}} \leq q_3 + 1,$ *and*

*(iii)*    $3 \cdot 2^{r_{3,o}} \leq q_3 + 1,$

*where $2^{r_{2,o}}, 2^{r_{q,o}},$ and, respectively, $2^{r_{3,o}},$ denote the largest powers of two less than or equal to the three quantities $2t, u,$ and, respectively, $v,$ with $2t = q_1 - q_2 - 2, u = u(q)$ the number of parts, among $q_4, q_5, \ldots q_r,$ that are equal to an odd integer $q,$ for any odd integer $q \neq 3,$ and, respectively, $v$ the number of parts, among $q_4, q_5, \ldots q_{r-1},$ that are equal to 3. (If any of these quantities is zero, the corresponding conditions (i), (ii), and (iii) do not apply.)*

*Then, for any $n \geq 6,$*

*$o'_e(n) = s_e(n)$ (see Definition 3.1).*

*(b) For any positive integer $n,$ $n \geq 6,$ let $\{o'_o(n)\}$ be the number of partitions of the form*

$$q_1 \geq q_2 = q_3 \geq q_4 \geq \ldots \geq q_r = 3$$

*with*

$$\sum_{j=1}^{r} q_j = n,$$

*and each part $q_j, j = 1, 2, \ldots, r$ odd, greater than or equal to 3, $r \geq 3$ (the partition must have at least three parts), the second largest part equal to the third, and the third largest part greater than the fourth largest unless the second, the third, and the fourth largest parts are all equal to 3, subject to the following conditions (merging caps):*

*(i)*      $2 \cdot 2^{r_{2,o}} \leq q_3 - 2,$

*(ii)*     $q \cdot 2^{r_{q,o}} \leq q_3 - 2,$ *and*

*(iii)*    $3 \cdot 2^{r_{3,o}} \leq q_3 - 2,$

*where $2^{r_{2,o}}, 2^{r_{q,o}},$ and, respectively, $2^{r_{3,o}},$ denote the largest powers of two less than or equal to the three quantities $2t, u,$ and, respectively, $v,$ with $2t = q_1 - q_2 - 2, u = u(q)$ the number of parts, among $q_4, q_5, \ldots q_r,$ that are equal to an odd integer $q,$ for any odd integer $q \neq 3,$ and, respectively, $v$ the number of parts, among $q_4, q_5, \ldots q_{r-1},$ that are equal to 3. (If any of these quantities is zero, the corresponding condition (i), (ii), or (iii) does not apply.)*

*Then, for any $n \geq 6,$*

*$o'_o(n) = s_o(n)$ (see Definition 3.1).*

*(c) For any positive integer $n,$ $n \geq 6,$*

$$o'_e(n) + o'_o(n) = s(n)$$

*(using the same notation of parts (a) and (b), and Definition 1.1).*

☐

(Proposition 3.3, here, corresponds to <span style="color:red">Proposition 3.4</span> in <span style="color:red">[H]</span>, p. 8.)

## 4   THE GENERALIZED PENTAGONAL NUMBER STRUCTURE OF THE BUTTERFLY SEQUENCE

The generalized pentagonal number structure of the sequence (1.10) is first observed as $s(n)$ is even *except when $n$ is a generalized pentagonal number, or a generalized pentagonal number plus two (or a generalized pentagonal number "with domino", for $n \geq 5$ and $n \neq 7$* (see Remark 1.2, here), and again displayed (equivalently, except for $n = 5$) as the first few terms of the butterfly subsequences $\{s_e(n)\}$ and $\{s_o(n)\}$ which add up to (1.10) (see (3.3) and (3.4) in the Definition 3.1, here), are equal for all $n \geq 6$, *except for*

$$n = 9, 12, 14, 15, 17, 22, 24, 28, 35, 37, 42, 51, \dots . \tag{4.1}$$

The special features of $n = 5$ and $n = 7$ are for the overlapping of the second difference of the sequence (1.10) (see (2.11), here).

*These exceptions, the inputs in (4.1) of the sequences (1.10), (3.3), and (3.4), are the generalized pentagonal numbers*, greater than 5, except, as noted, for $n = 7$, *together with the generalized pentagonal numbers plus two units (*what we call *the generalized pentagonal number plus domino)*, again larger than 5, *in* (3.3) *and* (3.4) (except for 7).

The following definitions (Definition 4.1, 4.2, and 4.4, below) of *pentagonal butterfly partitions, generalized pentagonal butterfly partitions, pentagonal butterfly partitions with domino, generalized pentagonal butterfly partitions with domino, and, respectively, non-pentagonal butterfly partitions,* together with the related theorem (Theorem 4.5, below) and corollary (Corollary 4.6, below) will show that these initial observations extend to all integers $n \geq 6$, $n \neq 7$. The proof of the Theorem 4.4, below, contains *non-pentagonal partition bijections* that show how the pentagonal numerical structure extends from the inputs of the partitions in (1.10), (3.1) and (3.2) to the Young configurations of the outputs.

**Definition 4.1:** For any integer $n \geq 6$ and $h \geq 3$, we define the *pentagonal partition* of $n$, or the *butterfly pentagonal partition* of $n$, to be the partition

$$p_1 = 2h - 1 > p_2 = 2h - 2 > p_3 = 2h - 3 > \cdots > p_{h-1} = h + 1 \ > p_h = h, \qquad (4.2)$$

with the pentagonal number sum

$$\sum_{i=1}^{h} p_i = \frac{3}{2} h^2 - \frac{1}{2} h = n$$

(see Figure 9, below), and, similarly, for any integer $n \geq 6$ and $h \geq 3$, we define the *generalized pentagonal, and not just pentagonal, partition* of $n$, or the *butterfly generalized pentagonal, and not just pentagonal, partition* of $n$, to be the partition

$$p_1 = 2h > p_2 = 2h - 1 > p_3 = 2h - 2 > \cdots > p_{h-1} = h + 2 \ > p_h = h + 1, \qquad (4.3)$$

with the generalized pentagonal number sum

$$\sum_{i=1}^{h} p_i = \frac{3}{2} h^2 + \frac{1}{2} h = n$$

(see Figure 10, below).

$\square$

(For the readers of [H], Definition 4.1, here, corresponds to part (a) of Definition 4.1 in [H], p. 8.)

**Definition 4.2:** For any integer $n \geq 6$ and $h \geq 3$, we define the *butterfly pentagonal partition with domino* of $n$ to be the partition

$$p_1 = 2h - 1 > p_2 = 2h - 2 > p_3 = 2h - 3 > \cdots > p_{h-1} = h + 1 \ > p_h = h >$$
$$> p_{h+1} = 2, \qquad (4.4)$$

with the sum equal to the following pentagonal number *plus two units*

$$\sum_{i=1}^{h+1} p_i = \frac{3}{2} h^2 - \frac{1}{2} h + 2 = n$$

(see Figure 11, below), and, similarly, for any integer $n \geq 6$ and $h \geq 2$, we define the *butterfly generalized pentagonal, and not just pentagonal, partition with domino* of $n$ to be the partition

$$p_1 = 2h > p_2 = 2h - 1 > p_3 = 2h - 2 > \cdots > p_{h-1} = h + 2 \ > p_h = h + 1 >$$

$$> p_{h+1} = 2, \tag{4.5}$$

with the sum equal to the following generalized pentagonal number *plus two units*

$$\sum_{i=1}^{h} p_i = \frac{3}{2} h^2 + \frac{1}{2} h + 2 = n$$

(see Figure 12, below).

☐

(For the readers of [H], Definition 4.2, here, does not correspond to Definition 4.2 in [H]. It is instead incorporated in part (b) of Definition 4.1 in [H], p. 8.)

**Remark 4.3:** Note that, while (4.2), (4.3) and (4.4) are defined for $h \geq 3$, the *generalized pentagonal, and not just pentagonal, butterfly partition with domino* (4.5) is defined for $h \geq 2$, which relates to the "exception of exceptions", $n \neq 7$, mentioned in the first paragraph of the present section.

☐

(For the readers of [H], Remark 4.3, here, does not appear in [H], and thus does not correspond to Remark 4.3 in [H]. Instead, Remark 4.3 in [H] corresponds to Remark 4.5, below.)

The following Young diagrams offer a visualization of the partitions (4.2), (4.3), (4.4), and (4.5) with the pentagonal structure (as a square-triangular partition structure) in black (the $h$ by $h$ squares) and blue (the $h-1$ by $h-1$ triangles, on the left, and the $h$ by $h$ triangles, on the right), and the dominoes in red.

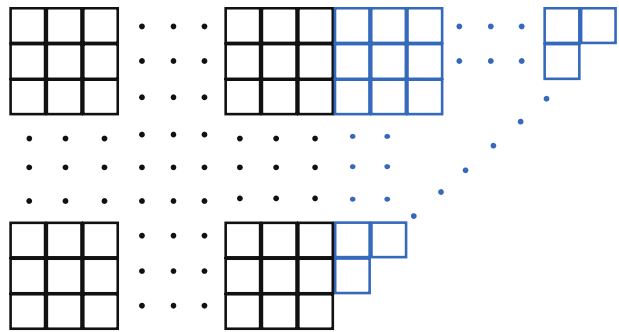

**Figure 9:** a butterfly pentagonal partition of $n = 12, 22, 35, 51, \dots$ .

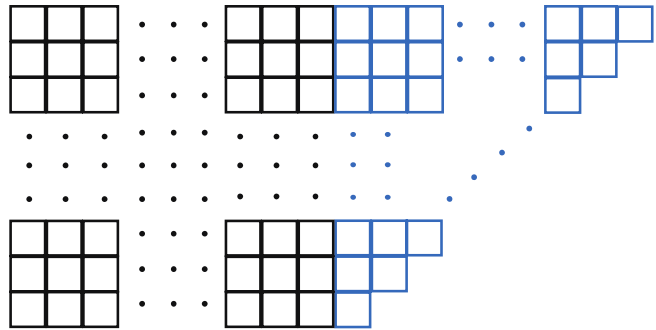

**Figure 10:** a generalized pentagonal butterfly partition (which is not a pentagonal butterfly partition) of $n = 15, 26, 40, 57, \ldots$.

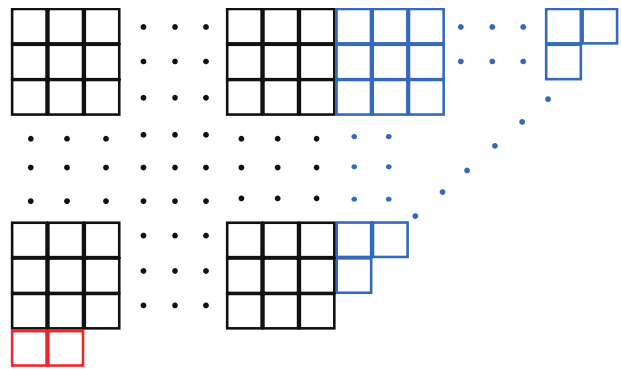

**Figure 11:** a pentagonal butterfly partition of $n = 14, 24, 37, 53, \ldots$ with the additional domino part in red.

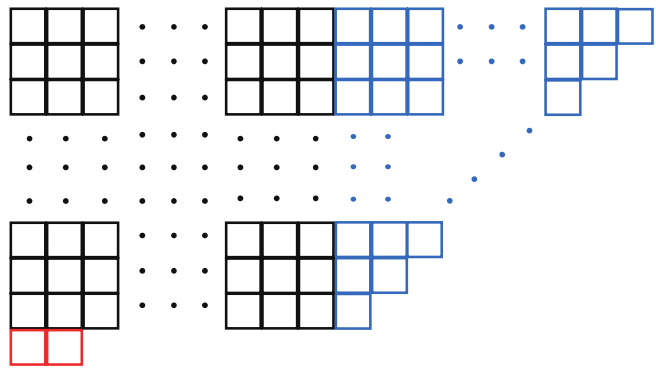

**Figure 12:** a generalized pentagonal butterfly partition (which is not a butterfly pentagonal partition) of $n = 9, 17, 28, 42, \ldots$ with the additional domino part for $n > 9$, in red.

**Definition 4.4:** (a) For any integer $n \geq 6$ and $h \geq 3$, let $A_e(n; h)$, and, respectively, $A_o(n; h)$, be the set of all the butterfly partitions of $n$ with the second part even, and, respectively, odd, *at least $h$ parts greater than or equal to $h + 1$, the largest $h$ parts* consecutive, and the smallest part, *other than 2,* equal to $h$. In detailed notation, the set $A_e(n; h)$ is composed of the partitions of $n$ of the form

$$p_1 = 2r + 1 > p_2 = 2r > p_3 = 2r - 1 > p_4 = 2r - 2 > \cdots > p_h = 2r - h + 2 > \cdots$$

$$\cdots > p_{k-1} > p_k = h, \tag{4.6}$$

with $h \geq 3$, together with the partitions of $n$ of the form

$$p_1 = 2r + 1 > p_2 = 2r > p_3 = 2r - 1 > p_4 = 2r - 2 > \cdots > p_h = 2r - h + 2 > \cdots$$

$$\cdots > p_{k-1} = h > p_k = 2. \tag{4.7}$$

Note that the set $A_e(n; h)$ is empty for $n < \frac{3}{2} h^2 + \frac{5}{2} h)$ (for the partition in $A_e(n; h)$ with the largest $h$ is $p_1 = 2h + 1 > p_2 = 2h > p_3 = 2h - 1 > \cdots > p_h = h + 2 > p_{h+1} = h)$, and note that, since $h \geq 3$ and there must be at least $h$ parts greater than or equal to $h + 1$, we have $r \geq 3$ in (4.6), $r \geq 4$ in (4.7), $n \geq 21$ in (4.6) (for the partition $7 > 6 > 5 > 3$ is the smallest such partition, for which the butterfly sequence output $s(21)$ is the third term equal to 2 in the sequence (1.10)), and $n \geq 23$ in (4.7) (for the partition $7 > 6 > 5 > 3 > 2$ is the smallest such partition, for which the butterfly sequence output $s(23)$ is the fourth term equal to 2 in the sequence (1.10)).

Similarly, in detailed notation, the set $A_o(n; h)$ is composed of the partitions of $n$ of the form

$$p_1 = 2r > p_2 = 2r - 1 > p_3 = 2r - 2 > p_4 = 2r - 3 > \cdots > p_h = 2r - h + 1 > \cdots$$

$$\cdots > p_{k-1} > p_k = h, \tag{4.8}$$

with $h \geq 3$, together with the partitions of $n$ of the form

$$p_1 = 2r > p_2 = 2r - 1 > p_3 = 2r - 2 > p_4 = 2r - 3 > \cdots > p_h = 2r - h + 1 > \cdots$$

$$\cdots > p_{k-1} = h > p_k = 2. \tag{4.9}$$

Note that the set $A_o(n; h)$ is empty for $n < \frac{3}{2} h(h + 1)$ (for the partition in $A_o(n; h)$ with the largest $h$ is $p_1 = 2h > p_2 = 2h - 1 > p_3 = 2h - 2 > \cdots > p_h = h + 1 > p_{h+1} = h)$, and note that, since h $\geq 3$ and there must be at least $h$ parts greater than or equal to $h + 1$, we have $r \geq 3$ in both (4.8) and (4.9), and $n \geq 18$ in (4.8) (for the partition $6 > 5 > 4 > 3$ is the smallest such partition, for which the butterfly sequence output $s(18)$ is the first term equal to 2 in the sequence (1.10)), and $n \geq 20$ in (4.9) (for the partition $6 > 5 > 4 > 3 > 2$ is the smallest such partition, for which the butterfly sequence output $s(20)$ is the second term equal to 2 in the sequence (1.10)).

(b) For any integer $n \geq 6$ and $h \geq 3$, let $B_e(n; h)$, and, respectively, $B_o(n; h)$, be the set of all the butterfly partitions of $n$ with the second part even, and, respectively, odd, at least $h$ parts, all parts, except the part 2 (if the 2 part is present), greater than $h$, the largest $h$ parts of the partition consecutive and greater than $h + 1$, and, if the number of parts of the partition is greater than $h$, the $(h + 1)$-th largest part at least two units less than the $h$-th part. In detailed notation, the set $B_e(n; h)$ is composed of the partitions of $n$ of the form

$$p_1 = 2r + 1 > p_2 = 2r > p_3 = 2r - 1 > \cdots > p_h = 2r - h + 2, \qquad (4.10)$$

the similar partitions with the additional part 2

$$p_1 = 2r + 1 > p_2 = 2r > p_3 = 2r - 1 > \cdots > p_h = 2r - h + 2 > p_{h+1} = 2, \qquad (4.11)$$

with $r > h - \frac{1}{2}$ in both (4.10) and (4.11) so that $2r - h + 2 > h + 1$, the partitions of $n$ of the form

$$p_1 = 2r + 1 > p_2 = 2r > p_3 = 2r - 1 > p_4 = 2r - 2 > \cdots > p_h = 2r - h + 2 > p_{h+1} > \cdots$$
$$\cdots > p_k \geq h + 1, \qquad (4.12)$$

together with the partitions of $n$ of the form

$$p_1 = 2r + 1 > p_2 = 2r > p_3 = 2r - 1 > p_4 = 2r - 2 > \cdots > p_h = 2r - h + 2 > p_{h+1} > \cdots$$
$$\cdots > p_{k-1} = h + 1 > p_k = 2, \qquad (4.13)$$

with $h \geq 3$, $r > h - \frac{1}{2}$ and $p_{h+1} < 2r - h + 1$ in both (4.12) and (4.13).

Note that the set $B_e(n; h)$ is empty for $n < \frac{3}{2}h(h + 1)$ (for the partition in $B_e(n; h)$ with the largest $h$ is $p_1 = 2h + 1 > p_2 = 2h > p_3 = 2h - 1 > \cdots > p_h = h + 2$) and note that, since $h \geq 3$, we have r $\geq 3$ in (4.10), (4.11), (4.12) and (4.13), $n \geq 18$ in (4.10) and (4.11) (for the partition $7 > 6 > 5$ is the smallest such partition, or, equivalently, for the butterfly sequence output $s(18)$ is the first term equal to 2 in the sequence (1.10)), and $n \geq 20$ in (4.12) and (4.13) (for the partition $7 > 6 > 5 > 2$ is the smallest such partition, or, equivalently, for the butterfly sequence output $s(20)$ is the second term equal to 2 in the sequence (1.10)).

Similarly, in detailed notation, the set $B_o(n; h)$ is composed of the partitions of $n$ of the form

$$p_1 = 2r > p_2 = 2r - 1 > p_3 = 2r - 2 > \cdots > p_h = 2r - h + 1, \qquad (4.14)$$

the similar partitions with the additional part 2

$$p_1 = 2r > p_2 = 2r - 1 > p_3 = 2r - 2 > \cdots > p_h = 2r - h + 1 > p_{h+1} = 2, \qquad (4.15)$$

with $r > h$ in both (4.14) and (4.15) so that $2r - h + 1 > h + 1$, the partitions of $n$ of the form

$$p_1 = 2r > p_2 = 2r - 1 > p_3 = 2r - 2 > \cdots > p_h = 2r - h + 1 > p_{h+1} > \cdots$$
$$\cdots > p_k \geq h + 1, \qquad (4.16)$$

together with the partitions of $n$ of the form

$$p_1 = 2r > p_2 = 2r - 1 > p_3 = 2r - 2 > \cdots > p_h = 2r - h + 1 > p_{h+1} > \cdots$$
$$\cdots > p_{k-1} = h + 1 > p_k = 2, \qquad (4.17)$$

with $h \geq 3$, $r > h$ and $p_{h+1} < 2r - h$ in (4.16) and (4.17).

Note that the set $B_o(n; h)$ is empty for $n < \frac{3}{2}h^2 + \frac{7}{2}h$ (for the partition in $B_o(n; h)$ with the largest $h$ is $p_1 = 2h + 2 > p_2 = 2h + 1 > p_3 = 2h > \cdots > p_h = h + 3$) and note that, since $h \geq 3$, we have also r $\geq 4$ in (4.14), (4.15), (4.16) and (4.17), $n \geq 21$ in (4.14) and (4.15) (for the partition $8 > 7 > 6$ is the smallest such partition, or, equivalently, for the butterfly sequence output s(21) is the third term equal to 2 in the sequence (1.10)), and $n \geq 23$ in (4.16) and (4.17) (for the partition $8 > 7 > 6 > 2$ is the smallest such partition, or, equivalently, for s(23) is the fourth term equal to 2 in the sequence (1.10)).

□

(For the readers of [H], Definition 4.4, here, corresponds to part to Definition 4.1 in [H], pp. 8 - 10.)

We call the partitions in the set $A_e(n; h) \cup A_o(n; h)$ the *non-pentagonal butterfly partitions of $n$ with the smallest part, other than the 2 part, equal $h$ and the $h$ largest parts consecutive and greater than $h$*, or, in short, as for the Young diagram visualizations of such partitions, the *non-pentagonal butterfly partitions of $n$ with the horizontal bar $h$*, and the partitions in the set $B_e(n; h) \cup B_o(n; h)$ the *non-pentagonal butterfly partitions of $n$ with the largest $h$ parts of the partition consecutive and greater than $h$, and with the smallest part, other than the 2 part, greater than or equal to $h + 1$*, or, in short, the *non-pentagonal butterfly partitions of $n$ with the vertical bar $h$* (see Figure 13, below). More specifically, we call the partitions in (4.6) the *non-pentagonal butterfly partitions of $n$ with the second part even and the horizontal bar $h$*, the partitions in (4.7) the *non-pentagonal butterfly partitions of $n$ with the second part even, the horizontal bar $h$, and a domino,* the partitions in (4.8) *non-pentagonal butterfly partitions of $n$ with the second part odd and the horizontal bar $h$*, and the partitions in (4.9) *non-pentagonal butterfly partitions of $n$ with the second part odd, the horizontal bar $h$, and a domino.*

(We use the words *non-pentagonal* to keep in mind that such non-pentagonal partitions do not have the Young diagram square-triangular configurations of pentagonal numbers.)

Similarly, we call the partitions in (4.10) and (4.12) the *non-pentagonal butterfly partitions of $n$ with the second part even and with the vertical bar $h$*, the partitions in (4.11) and (4.13) the *non-pentagonal butterfly partitions of $n$ with the second part even, the vertical bar $h$, and a domino,* the partitions in (4.14) and (4.16) *non-pentagonal butterfly partitions of $n$ with the second part odd and with the vertical bar $h$*, and the partitions in (4.15) and (4.17) the *non-pentagonal butterfly partitions of $n$ with the second part odd, the vertical bar $h$, and a domino.*

**Remark 4.5:** *Although there is a one-to-one correspondence between the set of generalized pentagonal numbers and the set of generalized butterfly pentagonal partitions, and a similar one-to-one correspondence between the set of generalized pentagonal numbers plus two units and the set of generalized butterfly pentagonal partitions with domino, the non-pentagonal butterfly partitions occur not only for integers $n$ that are not generalized pentagonal and are not generalized pentagonal plus two units, but also for all generalized pentagonal and generalized pentagonal plus two integers $n \geq 24$. For example, the pentagonal number plus two units $24$ has the non-pentagonal butterfly partition $8 > 7 > 6 > 3$ with horizontal bar 3, the corresponding non-pentagonal butterfly partition $9 > 8 > 7$ with vertical bar 3, and the pentagonal butterfly partition with domino $7 > 6 > 5 > 4 > 2$.*

□

(For the readers of [H], Remark 4.5, here, corresponds to Remark 4.3 in [H], p. 10.)

**Theorem 4.6:** *Using the same notation of Definition 3.1, Proposition 3.4 and Proposition 3.5, here, for any integer $n \geq 6$,*

   a) *the outputs of the butterfly sequence with the second largest part even are equal to the outputs of the butterfly sequence with the second largest part odd, that is $\mathsf{s_e}(n) = \mathsf{s_o}(n)$, for all $n \neq \frac{1}{2}(3t^2 + t + 4)$, $n \neq \frac{1}{2}(3(t+1)^2 - t - 1)$, $n \neq \frac{1}{2}(3(t+1)^2 - t + 3)$ and $n \neq \frac{1}{2}(3(t+1)^2 + t + 1)$,*

   b) *the outputs of the butterfly sequence with the second largest part even are one unit less than the outputs of the butterfly sequence with the second largest part odd, that is $\mathsf{s_e}(n) = \mathsf{s_o}(n) - 1$, for all $n = \frac{1}{2}(3t^2 + t + 4)$ and for all $n = \frac{1}{2}(3(t+1)^2 + t + 1)$, and*

c) *the outputs of the butterfly sequence with the second largest part even are one unit more than the outputs of the butterfly sequence with the second largest part odd, that is* $s_e(n) = s_o(n) + 1$, *for all* $n = \frac{1}{2}(3(t+1)^2 - t - 1)$ *and for all*

$n = \frac{1}{2}(3(t+1)^2 - t + 3)$,

*where $t$ is any integer greater than or equal to* 2.

*Proof.* For any integer $n \geq 6$, and for all $h \geq 3$, we start with $h = 3$, and we continue, term-by term, with $h = 4, 5, 6, 7, \dots$, (as far as the size of $n$ allows), setting two bijections, the first between the set $A_e(n; h)$ and the set $B_o(n; h)$, and the second bijection between the set $A_o(n; h)$ and the set $B_e(n; h)$. The first bijection, for $h = 3$, injects each partition of the form (4.6), and, respectively, (4.7), into a partition of the form (4.14) or (4.16), and, respectively, of the form (4.15) or (4.17), by removing the part equal to $h = 3$ (removing the horizontal bar $h = 3$) from (4.6), and, respectively, (4.7), and adding to the result one unit to each of the largest three (as $h = 3$) parts of the results, thus obtaining a partition of the form (4.14) or (4.16), and, respectively, (4.15) or (4.17). Vice versa, the first bijection removes one unit from each of the largest three (as $h = 3$) parts (it removes the vertical bar $h = 3$ in the partition on the right in Figure 13, below) in (4.14) or (4.16), and, respectively, in (4.15) and (4.17), and adds to the result a part equal to $h = 3$, thus reconstructing the original partitions. Similarly, the second bijection, again for $h = 3$, injects each partition of the form (4.8), and, respectively, (4.9), into a partition of the form (4.10) or (4.12), and, respectively, (4.11) or (4.13), by removing the part equal to $h = 3$ (the second bijection removes the horizontal bar $h = 3$) from (4.8), and, respectively, (4.9), and adding one unit to each of the largest three (as $h = 3$) parts of the results, thus obtaining a partition of the form (4.10) or (4.12), and, respectively, (4.11) or (4.13). Vice versa, the second bijection removes one unit from each of the largest $h$ parts (removes the vertical bar $h = 3$ in the partition on the right in Figure 13, below) in (4.14) or (4.16), and respectively, (4.15) and (4.17), and adds a part equal to $h = 3$ to the result, thus recovering the original partitions. After the two bijections are set for the partitions $h = 3$, note that the remaining non-pentagonal partitions in $A_e(n; h) \cup A_o(n; h)$ have the smallest part, other than 2, greater than or equal to 4, and note that the remaining non-pentagonal partitions in $B_e(n; h) \cup B_o(n; h)$ have the four largest parts consecutive, so we can set the two bijections for $h = 4$, step-by-step, as in the above ($h = 3$) procedure, except with $h = 4$ in place of $h = 3$. The procedure is then repeated for with $h = 5, 6, 7, \dots$, as far as the size of $n$ requires ($n < \frac{3}{2}h^2 + \frac{5}{2}h$ for the partitions in (4.6) and (4.7), and $n < \frac{3}{2}h(h+1)$ for the partitions in (4.8) and (4.9)). Since each butterfly partition (whether with the second part even, or with the second part odd), is necessarily pentagonal, generalized pentagonal, pentagonal with domino, generalized pentagonal with domino, or non-pentagonal, the results of the theorem are then a straightforward consequence of the cyclical order of appearance of the butterfly partitions as $n \geq 9$ increases (the four cycle a) generalized pentagonal with domino, b) pentagonal, c) pentagonal with domino, d) generalized pentagonal etc., below), combined with the 2-cycle of the second largest part of the corresponding partition ("odd, even, even, odd, etc. …" starting with the partition $4 > 3 > 2$ for $n = 9$). We do indeed identify one pentagonal, generalized (and not just pentagonal) pentagonal, pentagonal with domino, or generalized (and not just pentagonal) partition with domino for each pentagonal, generalized pentagonal, pentagonal plus two, or generalized pentagonal plus two number in the following order: a) $n = 9 = \frac{1}{2}(3t^2 + t + 4)$, a generalized pentagonal number plus two, with t = 2, for the butterfly generalized pentagonal partition with domino $4 > 3 > 2$; b) $n = 12 = \frac{1}{2}(3(t + 1)^2 - t - 1)$, the third pentagonal number, with t = 2, for the pentagonal butterfly partition $5 > 4 > 3$; c) $n = 14 = \frac{1}{2}(3(t + 1)^2 - t - 1) + 2 = \frac{1}{2}(3(t + 1)^2 - t + 3)$, the third pentagonal number plus two, with t = 2, for the butterfly pentagonal partition with domino $5 > 4 > 3 > 2$; d) $n = 15 = \frac{1}{2}(3(t + 1)^2 + t + 1)$, the third generalized (and not just pentagonal) pentagonal number, with t = 2, for the butterfly generalized pentagonal partition $6 > 5 > 4$; etcetera with t = $3, 4, 5, 6, 7, … .$

$\square$

(For the readers of [H], Theorem 4.6, here, above, corresponds to Theorem 4.4 in [H], p. 10. On the other hand, Figure  13, here, below, corresponds to Figure 2 in [H], p. 10. Figures $7 – 12$, here, above, do not appear in the compact version, [H], of the present paper.)

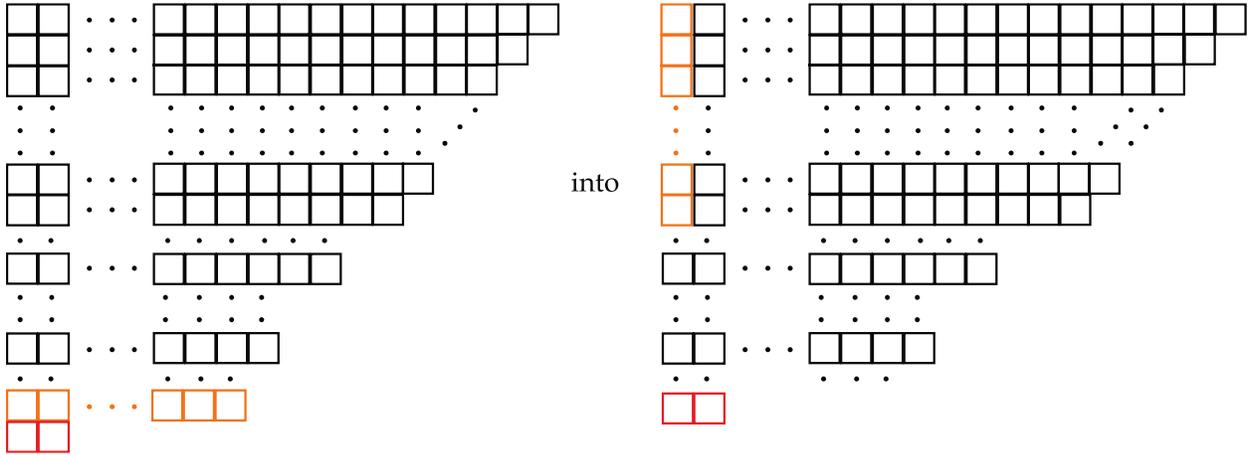

**Figure 13:** a non-pentagonal butterfly partition with the horizontal part *h* in orange, and its corresponding non-pentagonal butterfly partition with the vertical bar *h* also in orange, both with (red) dominoes. (In both partitions on the left and on the right, the second part may be even or odd.) Similar configurations occur without the red dominoes.

**Corollary 4.7:** *Using the same notation of Definition 3.1, Proposition 3.4 and Proposition 3.5, for any integer $n \geq 6$,*

a) *the outputs of the butterfly sequence are equal to $\mathsf{s}(n) = 2\mathsf{s}_\mathsf{e}(n)$ and, as such, are even for all $n \neq \frac{1}{2}(3t^2 + t + 4)$, $n \neq \frac{1}{2}(3(t+1)^2 - t - 1)$, $n \neq \frac{1}{2}(3(t+1)^2 - t + 3)$ and $n \neq \frac{1}{2}(3(t+1)^2 + t + 1)$,*

b) *the outputs of the butterfly sequence are equal to $\mathsf{s}(n) = 2\mathsf{s}_\mathsf{e}(n) + 1$ and, as such, are odd for all $n = \frac{1}{2}(3t^2 + t + 4)$ and for all $n = \frac{1}{2}(3(t+1)^2 + t + 1)$, and*

c) *the outputs of the butterfly sequence are equal to $\mathsf{s}(n) = 2\mathsf{s}_\mathsf{e}(n) - 1$ and, as such, are odd for all $n = \frac{1}{2}(3(t+1)^2 - t - 1)$ and for all $n = \frac{1}{2}(3(t+1)^2 - t + 3)$,*

*where t is any integer greater than or equal to 2.*

□

(For the readers of [H], Corollary 4.7, here, corresponds to Corollary 4.5 in [H], p. 11.)

Applying Corollary 1.4 to Theorem 4.6, here, we also have the following corollary.

**Corollary 4.8:** *For $n \geq 6$, setting $\mathsf{e}(n)$, respectively, $\mathsf{o}(n)$, the number of partitions of n with parts $p_i$, $i = 1, 2, \dots, k$, such that (a) $k \geq 3$, (b) the first three largest parts are even, respectively, odd, and equal, (c) the remaining parts, if present, are distinct, (d) the fourth largest part, if* *present, is two units smaller than the third largest part, and (e) the smallest part is not 1, we see that*

a) $e(n) = o(n)$, *for all* $n \neq \frac{1}{2}(3t^2 + t + 4)$, $n \neq \frac{1}{2}(3(t+1)^2 - t - 1)$, $n \neq \frac{1}{2}(3(t+1)^2 - t + 3)$ *and* $n \neq \frac{1}{2}(3(t+1)^2 + t + 1)$,

b) $e(n) = o(n) - 1$, *for all* $n = \frac{1}{2}(3t^2 + t + 4)$ *and for all* $n = \frac{1}{2}(3(t+1)^2 + t + 1)$, *and*

c) $e(n) = o(n) + 1$, *for all* $n = \frac{1}{2}(3(t+1)^2 - t - 1)$ *and for all*

   $n = \frac{1}{2}(3(t+1)^2 - t + 3)$,

*where t is any integer greater than or equal to* 2.

□

(For the readers of [H], Corollary 4.8, here, corresponds to Corollary 4.6 in [H], p. 11.)

Setting, for $n \geq 6$, $e'(n)$, the number of partitions of the form (1.14) (see Corollary (1.5) *with an even number of parts* ($h$ even in Corollary (1.5)), and, respectively, setting $o'(n)$ the number of partitions of the form (1.14) *with an odd number of parts* ($h$ odd), we have the following statement, similar to the thesis of Corollary 4.8 (except for the signs in b) and c) for $e'(n) = o(n)$, and $o'(n) = e(n)$).

**Corollary 4.9:** *For* $n \geq 6$,

a) $e'(n) = o'(n)$, *for all* $n \neq \frac{1}{2}(3t^2 + t + 4)$, $n \neq \frac{1}{2}(3(t+1)^2 - t - 1)$, $n \neq \frac{1}{2}(3(t+1)^2 - t + 3)$ *and* $n \neq \frac{1}{2}(3(t+1)^2 + t + 1)$,

b) $e'(n) = o'(n) + 1$, *for all* $n = \frac{1}{2}(3t^2 + t + 4)$ *and for all* $n = \frac{1}{2}(3(t+1)^2 + t + 1)$, *and*

c) $e'(n) = o'(n) - 1$, *for all* $n = \frac{1}{2}(3(t+1)^2 - t - 1)$ *and for all*

   $n = \frac{1}{2}(3(t+1)^2 - t + 3)$,

*where t is any integer greater than or equal to* 2.

□

(For the readers of [H], Corollary 4.9, here, corresponds to Corollary 4.7 in [H], p. 11.)

Finally, setting, for $n \geq 6$, e"$(n)$, the number of partitions of the form (1.15) (see Corollary (1.6) *with an even number of parts* ($h$ even in Corollary (1.6)), and, respectively, setting o"$(n)$ the number of partitions of the form (1.15) *with an odd number of parts* ($h$ odd), we have the following corollary.

**Corollary 4.10:** *For $n \geq 6$,*

a) e"$(n)$ = o"$(n)$, *for all $n \neq \frac{1}{2}(3t^2 + t + 4)$,  $n \neq \frac{1}{2}(3(t+1)^2 - t - 1)$,  $n \neq \frac{1}{2}(3(t+1)^2 - t + 3)$ and $n \neq \frac{1}{2}(3(t+1)^2 + t + 1)$,*

b) e"$(n)$=o"$(n)$ − 1, *for all $n = \frac{1}{2}(3t^2 + t + 4)$ and for all $n = \frac{1}{2}(3(t+1)^2 + t + 1)$, and*

c) e"$(n)$=o"$(n)$ + 1, *for all $n = \frac{1}{2}(3(t+1)^2 - t - 1)$ and for all*

$n = \frac{1}{2}(3(t+1)^2 - t + 3)$,

*where $t$ is any integer greater than or equal to 2.*

□

(For the readers of [H], Corollary 4.10, here, corresponds to Corollary 4.8 in [H], p. 11.)

# 5   SPECIALIZATIONS OF THE JACOBI TRIPLE PRODUCT, AND ALGORITHMS TO COMPUTE THE BUTTERFLY SEQUENCE AND RELATED SEQUENCES IN TERMS OF PENTAGONAL SEQUENCES AND TRIANGULAR NUMBERS

We begin by applying Euler's Pentagonal Number Theorem (cf. [A2], Theorem 1.6 and Corollaries 1.7 and 1.8) to the generating function of the sequence of partitions with distinct parts (2.1):

$$\sum_{m=0}^{\infty} q(m)x^m = \prod_{n=1}^{\infty}(1+x^n) = \prod_{n=1}^{\infty}\frac{1}{(1-x^n)}\,(1-x^{2n}) =$$

$$= \left(\sum_{h=0}^{\infty} p(h)x^h\right)\left(1 + \sum_{k=1}^{\infty}(-1)^k[x^{3k^2-k} + x^{3k^2+k}]\right) = \qquad (5.1)$$

$$= \sum_{h=0}^{\infty} p(h)x^h + \sum_{h=0}^{\infty}\sum_{k=1}^{\infty}(-1)^k[p(h)x^{h+3k^2-k} + p(h)x^{h+3k^2+k}]$$

with $q(m)$, $m = 0, 1, 2, \ldots$ , the number of partitions with distinct parts of $m$ (as in (1.1)), and $p(h)$, $h = 0, 1, 2, \ldots$ , the number of partitions of $h$, we obtain

$$q(m) = p(m) + \sum_{k=1}^{\infty} (-1)^k [\, p(m - 3k^2 + k) + p(m - 3k^2 - k) \,], \qquad (5.2)$$

setting $p(m) = 0$ for $m < 0$. Similar identities are easily obtained for the first and second (respectively) difference sequences of the sequence of partitions with distinct parts, by multiplying (5.1) by $1 - x$, and, respectively, by $1 - 2x + x^2$, and incorporating the two respective factors into the partition generating function factor of (5.1). Using the same notation (1.2) and (1.3) for the sequence of partitions with distinct parts and the two largest parts consecutive (1.2) and for the butterfly sequence (1.10), using the differential notation $dp(h)$, $h = 0, 1, 2, \ldots$ , for the number of partitions of $h$ that do not contain 1 as a part (see, for example, sequence A002865 in the O. E. I. S.), and denoting $d^2 p(h)$, $h = 0, 1, 2, \ldots$ , the number of partitions of $h$ that do not contain 1 as a part and *with the largest part occurring at least twice* (sequence A053445 in the O. E. I. S.), we see that

$$r(m) = dp(m) + \sum_{k=1}^{\infty} (-1)^k [\, dp(m - 3k^2 + k) + dp(m - 3k^2 - k) \,], \qquad (5.3)$$

setting $dp(h) = 0$ for $h < 0$, and

$$s(m) = d^2 p(m) + \sum_{k=1}^{\infty} (-1)^k [\, d^2 p(m - 2P_h) + d^2 p(m - 2P_h - 2h) \,], \qquad (5.4)$$

setting $d^2 p(h) = 0$ for $h < 0$.

On the other hand, incorporating the factors $1 - x$ and, respectively, $1 - 2x + x^2$ into the pentagonal number factor in (5.1) (multiplied left and right by the factor $1 - x$, and, respectively, $1 - 2x + x^2$, leads to the following series identities

$$\sum_{m=0}^{\infty} r(m) x^m =$$

$$= \left( \sum_{h=0}^{\infty} p(h) x^h \right) (1 + \sum_{k=1}^{\infty} (-1)^k [x^{3k^2 - k} + x^{3k^2 + k}] - x +$$

$$+ \sum_{k=1}^{\infty} (-1)^{k+1} [x^{3k^2 - k + 1} + x^{3k^2 + k + 1}]) = \qquad (5.5)$$

$$= \sum_{h=0}^{\infty} p(h)x^h + \sum_{h=0}^{\infty} \sum_{k=1}^{\infty} (-1)^k [p(h)x^{h+3k^2-k} + p(h)x^{h+3k^2+k}] - \sum_{h=0}^{\infty} p(h)x^{h+1} +$$

$$+ \sum_{h=0}^{\infty} \sum_{k=1}^{\infty} (-1)^{k+1} [p(h)x^{h+3k^2-k+1} + p(h)x^{h+3k^2+k+1}],$$

and, respectively,

$$\sum_{m=0}^{\infty} s(m)x^m = \left( \sum_{h=0}^{\infty} p(h)x^h \right) \cdot$$

$$\cdot (1 + \sum_{k=1}^{\infty} (-1)^k [x^{3k^2-k} + x^{3k^2+k}] - 2x + \sum_{k=1}^{\infty} 2(-1)^{k+1}[x^{3k^2-k+1} + x^{3k^2+k+1}] +$$

$$+ x^2 + \sum_{k=1}^{\infty} (-1)^k [x^{3k^2-k+2} + x^{3k^2+k+2}]) = \tag{5.6}$$

$$= \sum_{h=0}^{\infty} p(h)x^h + \sum_{h=0}^{\infty} \sum_{k=1}^{\infty} (-1)^k [p(h)x^{h+3k^2-k} + p(h)x^{h+3k^2+k}] - 2\sum_{h=0}^{\infty} p(h)x^{h+1} +$$

$$+ \sum_{h=0}^{\infty} \sum_{k=1}^{\infty} 2(-1)^{k+1}[p(h)x^{h+3k^2-k+1} + p(h)x^{h+3k^2+k+1}] + \sum_{h=0}^{\infty} p(h)x^{h+2} +$$

$$+ \sum_{h=0}^{\infty} \sum_{k=1}^{\infty} (-1)^k [p(h)x^{h+3k^2-k+2} + p(h)x^{h+3k^2+k+2}],$$

which, in turn, give

$$r(m) = p(m) - p(m-1) + \sum_{k=1}^{\infty} (-1)^k \{p(m-3k^2+k) + p(m-3k^2-k) +$$

$$-p(m-3k^2+k-1) + p(m-3k^2-k-1) \}, \tag{5.7}$$

and, respectively,

$$s(m) = p(m) - 2p(m-1) + p(m-2) +$$

$$+ \sum_{k=1}^{\infty} (-1)^k \{p(m-3k^2+k) + p(m-3k^2-k) - 2p(m-3k^2+k-1) + \tag{5.8}$$

$$-2p(m-3k^2-k-1) + p(m-3k^2+k-2) + p(m-3k^2-k-2) \}.$$

Besides keeping in mind that Euler's Pentagonal Number Theorem is a special case of the Jacobi triple product

$$\prod_{m=1}^{\infty} (1 - z^{2m}) \, (1 + z^{2m-1} y^2) \, (1 + z^{2m-1}/y^2) = \sum_{n=-\infty}^{\infty} z^{n^2} y^{2n} \tag{5.9}$$

(cf. [J], pp. 89-90, and [A1], Section 7.2), we can obtain a slightly different specialization of the Jacobi triple product (5.9) in the following two steps: first, set $y = z^{1/2}$ in (5.9), divide by 2 both sides, and multiply the first two sets of factors of the resulting triple product, thus writing

$$\prod_{m=1}^{\infty} (1 - z^{4m}) \, (1 + z^{2m}) = \frac{1}{2} \sum_{n=-\infty}^{\infty} z^{n^2+n} = \tag{5.10}$$

$$= \frac{1}{2} \{ 1 + \sum_{n=1}^{\infty} z^{n^2+n} + \sum_{n=1}^{\infty} z^{n^2-n} \} = \sum_{k=0}^{\infty} z^{k^2+k},$$

and, for the second step, set $z = x^{1/2}$ in (5.10) to write down a double product identity

$$\prod_{m=1}^{\infty} (1 + x^m) \, (1 - x^{2m}) = \sum_{k=0}^{\infty} x^{\frac{1}{2}k^2 + \frac{1}{2}k}. \tag{5.11}$$

Dividing both sides of (5.11) by

$$\prod_{m=1}^{\infty} (1 - x^{2m}),$$

and using the series expansion of what we call the *domino partition function* (as usual for the Young diagram visualization), we obtain

$$\sum_{m=0}^{\infty} q(m) x^m = \prod_{m=1}^{\infty} (1 + x^m) = \left( \prod_{m=1}^{\infty} \frac{1}{(1 - x^{2m})} \right) \left( \sum_{k=0}^{\infty} x^{\frac{1}{2}k^2 + \frac{1}{2}k} \right) = \tag{5.12}$$

$$= \left( \sum_{h=0}^{\infty} p(h) x^{2h} \right) \left( \sum_{k=0}^{\infty} x^{\frac{1}{2}k^2 + \frac{1}{2}k} \right) = \sum_{h=0}^{\infty} \sum_{k=0}^{\infty} p(h) x^{2h + \frac{1}{2}k^2 + \frac{1}{2}k}$$

and the component identities

$$q(m) = \sum_{k=0}^{\infty} p \left( \frac{1}{2} m - \frac{1}{4} k^2 - \frac{1}{4} k \right), \tag{5.13}$$

where, as usual, we set $p(h) = 0$, for all rational numbers $h$ except the nonnegative integers. Multiplying (5.12) by $1 - x$ and, respectively, $1 - 2x + x^2$, and incorporating the two respective factors into the generating function factor of (4.12), using the same notation as in (5.3) and (5.4) above, and setting $dp(h) = d^2p(h) = 0$, for all rational numbers h except the nonnegative integers, we see that

$$r(m) = \sum_{k=0}^{\infty} dp\left(\frac{1}{2}m - \frac{1}{4}k^2 - \frac{1}{4}k\right), \tag{5.14}$$

and

$$s(m) = \sum_{k=0}^{\infty} d^2p(h)\left(\frac{1}{2}m - \frac{1}{4}k^2 - \frac{1}{4}k\right) \tag{5.15}$$

As in (5.5) and (5.6), the incorporation of the factors $1 - x$ and, respectively, $1 - 2x + x^2$ into the triangular number factor of (5.12) (again, multiplying (5.12) by the factors $1 - x$ and, respectively, $1 - 2x + x^2$) leads to the series identities

$$\sum_{m=0}^{\infty} r(m)x^m = \left(\sum_{h=0}^{\infty} p(h)x^{2h}\right)\left(\sum_{k=0}^{\infty} x^{\frac{1}{2}k^2 + \frac{1}{2}k} - \sum_{k=0}^{\infty} x^{\frac{1}{2}k^2 + \frac{1}{2}k + 1}\right) = \tag{5.16}$$

$$= \sum_{h=0}^{\infty}\sum_{k=0}^{\infty}[p(h)x^{2h + \frac{1}{2}k^2 + \frac{1}{2}k} - p(h)x^{2h + \frac{1}{2}k^2 + \frac{1}{2}k + 1}],$$

and, respectively,

$$\sum_{m=0}^{\infty} s(m)x^m = \left(\sum_{h=0}^{\infty} p(h)x^{2h}\right)\left(\sum_{k=0}^{\infty} x^{\frac{1}{2}k^2 + \frac{1}{2}k} - \sum_{k=0}^{\infty} 2x^{\frac{1}{2}k^2 + \frac{1}{2}k + 1} + \sum_{k=0}^{\infty} x^{\frac{1}{2}k^2 + \frac{1}{2}k + 2}\right) = \tag{5.17}$$

$$= \sum_{h=0}^{\infty}\sum_{k=0}^{\infty}[p(h)x^{2h + \frac{1}{2}k^2 + \frac{1}{2}k} - 2p(h)x^{2h + \frac{1}{2}k^2 + \frac{1}{2}k + 1} + p(h)x^{2h + \frac{1}{2}k^2 + \frac{1}{2}k + 2}].$$

Taking the coefficients of $x^m$ in (5.16) and, respectively, in (5.17), we see that

$$r(m) = \sum_{k=0}^{\infty}[p\left(\frac{1}{2}m - \frac{1}{4}k^2 - \frac{1}{4}k\right) - p\left(\frac{1}{2}m - \frac{1}{4}k^2 - \frac{1}{4}k - 1\right)], \tag{5.18}$$

and, respectively,

$$s(m) = \sum_{k=0}^{\infty} [p\left(\frac{1}{2}m - \frac{1}{4}k^2 - \frac{1}{4}k\right) - 2p\left(\frac{1}{2}m - \frac{1}{4}k^2 - \frac{1}{4}k - 1\right) +$$

$$+ p\left(\frac{1}{2}m - \frac{1}{4}k^2 - \frac{1}{4}k - 2\right)]. \tag{5.19}$$

The identity (5.11) also gives

$$\prod_{m=1}^{\infty} (1 + x^m)(1 - x^{2m}) = \prod_{m=1}^{\infty} \frac{1}{(1 - x^{2m-1})}(1 - x^{2m}) =$$

$$= \left(\sum_{h=0}^{\infty} q(h)x^h\right)\left(1 + \sum_{k=1}^{\infty} (-1)^k [x^{3k^2-k} + x^{3k^2+k}]\right) = \tag{5.20}$$

$$= \sum_{h=0}^{\infty} q(h)x^h + \sum_{h=0}^{\infty}\sum_{k=1}^{\infty} (-1)^k q(h)[x^{h+3k^2-k} + x^{h+3k^2+k}] = \sum_{l=0}^{\infty} x^{\frac{1}{2}l^2 + \frac{1}{2}l},$$

from which we see a straightforward recursively algorithm *with pentagonal number sequences in the inputs of* $q(m)$ *and triangular numbers in the output*, as the expression

$$q(m) + \sum_{k=1}^{\infty} (-1)^k [q(m - 3k^2 + k) + q(m - 3k^2 - k)]$$

(in which, as usual, $q(h) = 0$ for $h$ negative) is equal to 1 if $m$ is a triangular number, and 0 otherwise. This algorithm extends to $r(m)$, $s(m)$, and beyond, by means of cyclotomic polynomials. Multiplying (5.20) by $1 - x$, we obtain

$$\sum_{h=0}^{\infty} r(h)x^h + \sum_{h=0}^{\infty}\sum_{k=1}^{\infty} (-1)^k r(h)[x^{h+3k^2-k} + x^{h+3k^2+k}] =$$

$$= \sum_{l=0}^{\infty} x^{\frac{1}{2}l^2 + \frac{1}{2}l} - \sum_{l=0}^{\infty} x^{\frac{1}{2}l^2 + \frac{1}{2}l + 1} = \tag{5.21}$$

$$= 1 + \sum_{l=2}^{\infty} [x^{\frac{1}{2}l^2 + \frac{1}{2}l} - x^{\frac{1}{2}l^2 + \frac{1}{2}l + 1}],$$

from which we see a similar *pentagonal and triangular* algorithm to compute $r(m)$, as the expression

$$r(m) + \sum_{k=1}^{\infty} (-1)^k [r(m - 3k^2 + k) + r(m - 3k^2 - k)]$$

(in which, as usual, $r(h) = 0$ for $h$ negative) is equal to 1, if $m = \frac{1}{2}l^2 + \frac{1}{2}l$, for $l = 0$ and $l \geq 2$, is equal to - 1, if $m = \frac{1}{2}l^2 + \frac{1}{2}l + 1$, for $l \geq 2$, and is equal to 0, if $m = 1$, or if $m$ or $m - 1$ are not triangular $\frac{1}{2}l^2 + \frac{1}{2}l$, for $l \geq 2$ (if $m \neq \frac{1}{2}l^2 + \frac{1}{2}l$ or $\frac{1}{2}l^2 + \frac{1}{2}l + 1$, for $l \geq 2$).

Again, multiplying (5.21) by $1 - x$, we obtain

$$\sum_{h=0}^{\infty} s(h)x^h + \sum_{h=0}^{\infty} \sum_{k=1}^{\infty} (-1)^k s(h)[x^{h+3k^2-k} + x^{h+3k^2+k}] =$$
$$= \sum_{l=0}^{\infty} \left[ x^{\frac{1}{2}l^2 + \frac{1}{2}l} - 2x^{\frac{1}{2}l^2 + \frac{1}{2}l + 1} + x^{\frac{1}{2}l^2 + \frac{1}{2}l + 2} \right] = \quad (5.22)$$
$$= 1 - x - x^2 + 2x^3 - 2x^4 + x^5 + \sum_{l=3}^{\infty} [x^{\frac{1}{2}l^2 + \frac{1}{2}l} - 2x^{\frac{1}{2}l^2 + \frac{1}{2}l + 1} + x^{\frac{1}{2}l^2 + \frac{1}{2}l + 2}],$$

from which we get a *pentagonal and triangular* algorithm to compute the butterfly sequence, as the expression

$$s(m) + \sum_{k=1}^{\infty} (-1)^k [s(m - 3k^2 + k) + s(m - 3k^2 - k)] \quad (5.23)$$

(in which, as usual, $s(h) = 0$ for $h$ negative) satisfies the following conditions:

(5.23) is equal to 1, if $m = \frac{1}{2}l^2 + \frac{1}{2}l$ for $l = 0$ and $l \geq 3$, or if $m = \frac{1}{2}l^2 + \frac{1}{2}l + 2$ for $l \geq 2$;

(5.23) is equal to $-1$, if $m = 1$ or $m = 2$;

(5.23) is equal to $-2$, if $m = \frac{1}{2}l^2 + \frac{1}{2}l + 1$ for $l \geq 2$;

and (5.24) is equal to 0, if $m \geq 6$ and $m, m - 1$ or $m - 2$ are not triangular for $l \geq 3$ ($m \neq \frac{1}{2}l^2 + \frac{1}{2}l$, $\frac{1}{2}l^2 + \frac{1}{2}l + 1$ or $\frac{1}{2}l^2 + \frac{1}{2}l + 2$, for $l \geq 3$).

Further multiplication of (5.22) by $1 + x + x^2$ leads to the similar identity, below, and a related similar algorithm

$$\sum_{h=0}^{\infty} t(h)x^h + \sum_{h=0}^{\infty} \sum_{k=1}^{\infty} (-1)^k t(h)[x^{h+3k^2-k} + x^{h+3k^2+k}] =$$

$$= \sum_{l=0}^{\infty} [x^{\frac{1}{2}l^2+\frac{1}{2}l} - 2x^{\frac{1}{2}l^2+\frac{1}{2}l+1} + x^{\frac{1}{2}l^2+\frac{1}{2}l+2} + x^{\frac{1}{2}l^2+\frac{1}{2}l+1} - 2x^{\frac{1}{2}l^2+\frac{1}{2}l+2} +$$

$$+ x^{\frac{1}{2}l^2+\frac{1}{2}l+3} + x^{\frac{1}{2}l^2+\frac{1}{2}l+2} - 2x^{\frac{1}{2}l^2+\frac{1}{2}l+3} + x^{\frac{1}{2}l^2+\frac{1}{2}l+4}] = \qquad (5.24)$$

$$= 1 - x^2 - x^4 + x^5 - x^9 + 2x^{10} - 2x^{11} + x^{12} +$$

$$+ \sum_{l=5}^{\infty} [x^{\frac{1}{2}l^2+\frac{1}{2}l} - x^{\frac{1}{2}l^2+\frac{1}{2}l+1} - x^{\frac{1}{2}l^2+\frac{1}{2}l+3} + x^{\frac{1}{2}l^2+\frac{1}{2}l+4}]$$

(see (2.9) and (2.10)), from which we get a *pentagonal and triangular* algorithm to compute the sequence of the number of partitions of n with odd parts larger or equal to 5, as the expression

$$t(m) + \sum_{k=1}^{\infty} (-1)^k [t(m-3k^2+k) + t(m-3k^2-k)] \qquad (5.25)$$

(in which, as usual, $t(h) = 0$ for $h$ negative) satisfies the following conditions:

(5.25) is equal to 1, if $m = \frac{1}{2}l^2 + \frac{1}{2}l$ for $l = 0$ and $l \geq 5$, if $m = \frac{1}{2}l^2 + \frac{1}{2}l + 2$ for $l = 0, 2, 4$, or if $m = \frac{1}{2}l^2 + \frac{1}{2}l + 4$ for $l \geq 5$;

(5.25) is equal to $-1$, if $m = \frac{1}{2}l^2 + \frac{1}{2}l + 1$ for $l = 1, 2$, or $l \geq 5$, or if $m = \frac{1}{2}l^2 + \frac{1}{2}l + 3$ for $l = 3$, or $l \geq 5$;

(5.25) is equal to 2, if $m = 10$;

(5.25) is equal to $-2$, if $m = 11$;

and (5.25) is equal to 0, for all the remaining nonnegative integers $m$.

In the expressions (5.21), (5.22), and (5.24), respectively, note the *wrinkled* polynomials 1, $1 - x - x^2 + 2x^3 - 2x^4 + x^5$, and $1 - x^2 - x^4 + x^5 - x^9 + 2x^{10} - 2x^{11} + x^{12}$, respectively, created by the multiplication of the series of triangular powers by the cyclotomic polynomials $1 - x$ and $1 + x + x^2$.